\documentclass[12pt,twoside]{article}
\usepackage{amsmath}
\usepackage{amssymb}
\usepackage{enumerate,epsfig}
\newtheorem{lemma}{Lemma}[section]
\newtheorem{theorem}{Theorem}[section]
\newtheorem{corollary}[lemma]{Corollary}

\newtheorem{prop}{Proposition}[section]

\newcommand{\R}{ {\mathbb R} }

\newcommand{\RR}{ {\mathbb R} }

\newcommand{\ep}{{\varepsilon}}
\newcommand{\pe}{{\varphi_\varepsilon}}
\newcommand{\fe}{{\bar\varphi_\varepsilon}}
\makeatletter
\@addtoreset{equation}{section}
\makeatother

\newcommand{\cqfd}{{\unskip\kern 6pt\penalty 500
\raise -2pt\hbox{\vrule\vbox to 6pt{\hrule width 6pt
\vfill\hrule}\vrule}\par}}
\newcommand{\eps}{{\varepsilon}}
\newcommand{\ind}{{\mathbb I}}

\begin{document}
\title{Small populations corrections for selection-mutation models} 

\date{}
\author{Pierre-Emmanuel Jabin$^{1}$}

\footnotetext[1]{CSCAMM and Department of Mathematics, University of
  Maryland, College Park, MD 20742-4015, USA \\
E-mail: \texttt{pjabin@cscamm.umd.edu}}

\maketitle{}

\begin{abstract}
  We consider  integro-differential models describing
  the evolution of a population structured 
  by a quantitative trait. Individuals interact competitively, creating a strong selection pressure on the population. On the other hand, mutations are assumed to be small. Following the formalism of~\cite{DJMP}, this creates
 concentration phenomena, typically consisting in a sum of Dirac masses
  slowly evolving in time. We propose a modification to those classical models that takes the effect of small populations into accounts and corrects some abnormal behaviours.

\end{abstract}

\bigskip


\section{Introduction}
\subsection{The first model}
We study  the dynamics of a population subject to mutations
and selection due to competition between individuals.
Each individual in the population is characterized by a
quantitative phenotypic trait $x$ 
(for example the size of the individual, their age at maturity, or their
rate of intake of nutrients). For simplicity, $x$ is taken here in $\R$ even though all the arguments could easily be extended to higher dimensional cases.

Probabilistic models are usually considered as the most realistic in that setting. They consist in life and death processes for each individual $X_i$, Poisson processes more precisely, with for instance birth rate $b(X_i)$ and a death rate which increases with the competition between individuals, for example
\[
d_i=d(X_i)+\sum_{j\neq i} I(X_i-X_j).
\]
When a birth occurs, it simply adds another individual with the same trait, except when a mutation takes place with small probability. In that case  the new individual has a different random trait, obtained through some distribution $K$. We refer to \cite{Me} and the references therein for a nice introduction to the probabilistic approach. 

Of course this is only one possible model and there are many variants. One can modify the interaction between individuals for instance by introducing explicit resources (with chemostat like interactions maybe). The competition could influence both the birth rate and the mortality rate...

When the total number of individuals is too large (it can easily reach $10^{10}-10^{12}$ for some micro-organisms), it becomes prohibitive to compute numerically the solution to this process. In that case one expects to be able to derive a deterministic model as a limit of large populations. Such of derivation was proved in \cite{CFM} and one obtains integro-differential equations like
\begin{equation}
\partial_t u(t,x)=\left(b(x)-d(x)-\int I(x-y)u(t,y)\,dy\right)\,u(t,x)
+M(u)(t,x), \label{ufirst} 
\end{equation}
where the mutation kernel is for instance
\[
M(f)(x)=\int_\R K(z)\,(b(x+z)\,f(x+z)-b(x)\,f(x))\,dz.
\]
\subsection{The scaling of fast reaction and small mutations}
Eq. \eqref{ufirst} is easy and fast to solve. However in most situations, small parameters appear and complicate that. This is a consequence of two small scales that are typical for this problem:
\begin{itemize}
\item The rate of mutation is very small in comparison to the reproduction rate. As we wish to see the evolution of traits generated by the mutations, this means that we need to rescale the equation in time and will hence get a large reproduction rate.
\item The size of each mutation is small.
\end{itemize} 
To simplify again our equations, let us now assume that $b=1$ and denote $r(x)=1-d(x)$ the reproduction rate of an individual without competition.

Taking the two scalings into account, one has in fact to deal with an equation like
\begin{equation}
\partial_t u_\eps(t,x)=\frac{1}{\ep}\left(r(x)-I\star u_\ep\right)\,u_\ep(t,x)
+M_\ep(u_\ep)(t,x), \label{u2} 
\end{equation}
where the mutation kernel now reads
\[
M_\ep(f)(x)=\frac{1}{\ep}\int_\R K(z)\,(b(x+\eps\,z)\,f(x+\eps z)-b(x)\,f(x))\,dz.
\]
Note that $M_\ep(f)$ is indeed of order $1$ if $b\,f\in C^1$ for instance.

Eq. \eqref{u2} is now much more delicate. The properties of the solution might depend on $\ep$ (its smoothness for instance) and solving it numerically can become again very costly if $\ep$ is too small (a typical value for many applications would be $\ep\sim 10^{-4}$). Therefore one would wish to derive a new model as $\ep\rightarrow 0$.

Eq. \eqref{u2} is strongly similar to a reaction-diffusion equation with a strong reaction term. However a crucial difference here is that the reaction term is non local. As we will see this completely changes the behaviour of the solution.

This scaling was introduced in \cite{DJMP} and formally studied there. We briefly reproduce the main argument here. The starting point is to introduce a large deviation scaling (which makes sense in view of the original probabilistic interpretation of the model)
\begin{equation}
\pe(t,x)=\ep\,\log u_\ep.\label{defphieps}
\end{equation}
Then it is easy to see that Eq. \eqref{u2} becomes
\begin{equation}
\partial_t \varphi_\eps=r-I\star u_\ep+H_\eps(\varphi_\eps),\label{phieps0}
\end{equation}
where
\begin{equation}
H_\eps(f)=\int_\R
K(z)\,\left(e^{(f(x+\eps\,z)-f(x))/\eps}-1\right)\,dz.\label{Heps0} 
\end{equation}
Using the theory of viscosity solutions to Hamilton-Jacobi equations, one can pass to the limit in \eqref{phieps0} and obtain for the limit $\varphi$ of $\pe$
\begin{equation}
\partial_t \varphi=r-I\star u+H(\partial_x\varphi),\label{phi0}
\end{equation}
with $u$ the weak limit (in the space of measure) of $u_\ep$ and 
\begin{equation}
H_\eps(f)=\int_\R
K(z)\,\left(e^{(f(x+\eps\,z)-f(x))/\eps}-1\right)\,dz.\label{H} 
\end{equation}
Of course at the limit, $\varphi$ and $u$ are no more connected by a relation like \eqref{defphieps}. Therefore Eq. \eqref{phi0} is no more closed and the question of how to recover $u$ from $\varphi$ is one of the main difficulty here, which is now fortunately better understood though.
\subsection{The problem with Eq. \eqref{u2}}
Let us focus here on one other delicate issue with this approach which is the main motivation for the current work.  In the scaling
under consideration, one has growth or decay of order
$\exp(C/\varepsilon)$. In particular one can see that the ratio between the maximal value of $u_\ep$ and the value at most other point is of this order  $\exp(C/\varepsilon)$. However if we come back to the starting point, which means a total population of $10^{10}-10^{12}$ and $\ep\sim 10^{-4}$, then there is an obvious problem. If (and it can be proved) $u_\ep$ is of order $\exp(C/\varepsilon)$ over a fixed interval of traits $I$ then the total population over this interval includes in fact much less than one individual!

This has several consequences. First of all it means that given the scales under consideration here, the limit of the probabilistic models to the deterministic one is not fully justified (the total number of individuals would have to be much higher with respect to $\ep$). This, in itself, could be ignored as anyway with $10^{12}$ individuals the obvious solution, use the probabilistic model, is not practical. Hence one could reasonably accept to still work with the deterministic model if its predictions were qualitatively in agreement with the behaviour of the stochastic one.  That is not the case.

The first indication of a serious flaw comes from the case where we do not put any mutations in the model. Making $M_\ep=0$ in \eqref{u2} still allows to perform the same analysis. In such a case, one would expect that nothing happens anymore at the limit as no evolution should be possible without mutations. However this is not the case and numerical simulations in particular show a remarquably similar behaviour between the cases with and without mutations (see \cite{PGau} for example). This phenomenon is entirely due to the persistence of very small subpopulations which should actually be extinct but are kept by our deterministic model and can therefore re-emerge later if the conditions are right.

A second important problem concerns the issue of branching. Biologically speaking branching is the process by which one population divides itself into two (and then possibly more) subpopulations. Mathematically one expects that at the limit, $u$ will be a sum of moving Dirac masses 
\[
u=\sum_i \alpha_i \delta_{x_i(t)}.
\]
Branching then corresponds to the case when one Dirac mass becomes two. This phenomenon is one of the main motivation to study models like \eqref{u2} instead of the adaptive dynamics approach for instance (see \cite{diekmann-04} for an introduction). At the limit \eqref{phi0}, branching occurs at infinite speed, {\em i.e.} if Dirac mass $\alpha \delta_{x(t)}$ divides itself after some time $t_0$ into $\alpha_1 \delta_{x_1(t)}+\alpha_2 \delta_{x_2(t)}$ then 
\[
\frac{d}{dt} |x_2(t)-x_1(t)|_{t=t0}=+\infty.
\] 
Instead probabilistic models predict a finite speed at branching...

Those qualitative disagreements are important and should be corrected. 
One would then hope to derive models that are both able of
dealing with very large populations and still treat correctly the
small subpopulations. It is for the moment a completely open question of how to keep the stochastic effects for the small populations. But one could at least try to
truncate the populations with less than 1 individual, which is expected to be enough to correct the qualitative flaws of the deterministic models. 

There are already some attempts in this direction, see \cite{PGau} and very recently \cite{MBPS}. The proposed correction in \cite{PGau} consists in studying equations like
\begin{equation}
\partial_t u_\eps(t,x)=\frac{1}{\ep}\left(r(x)-I\star u_\ep\right)\,u_\ep(t,x)
-\frac{1}{\ep}\,\sqrt{\frac{u_\ep}{\bar u_\ep}}+M_\ep(u_\ep)(t,x). \label{u2corr} 
\end{equation}
The added mortality term has the effect of killing (in times of order $\ep$) any population with a density less than $\bar u_\ep$. From the modeling point of view, this is quite satisfactory as it corrects most of the problems with \eqref{u2}.

Unfortunately the mathematical analysis of an equation like \eqref{u2corr} is for the moment untractable. The only situation that is understood
is when $\bar u_\ep$ is chosen like  $\exp(-\bar\varphi/\varepsilon)$. However this is exactly the scaling that was giving less than $1$ individual and it is not satisfactory. Instead given our scalings, one would like to work with $\bar u_\ep$ that are polynomial in $\ep$. 

However that would mean truncating every value less than $C\,\ep\log\ep$ in $\pe$, which implies at the limit, as $\varphi\leq 0$, everything...
%
\subsection{The proposed correction: Cooperative interactions}
%
We propose here a correction that is inspired from the case with sexual reproduction. In the present context of mostly asexual reproduction, it can however be better understood as taking into account some cooperative effects between the individuals. 

We can prove that it completely corrects all the abnormal behaviours of \eqref{u2} at the limit. In addition it is for the moment the only correction for which one can derive rigorously the limit. That means in particular that one can obtain numerical simulations for realistically low $\ep$.

There exists a well known phenomenon in the case of sexual reproduction which drives to extinction small populations. If the population is too small then the probability of meeting a partner is very low and hence the birth rate declines. This is the idea that we follow here.

Let us explain it first in the context of sexual reproduction. Consider a subpopulation with trait $x$, if the density of population for $x$, $u_\ep(x)$ is below a critical value $\bar u_\ep$ then we assume that the probability for an individual with trait $x$ to meet a partner is too low. Instead this individual will reproduce with an individual with a different trait $y$ such that $u_\ep(y)$ is large enough. Typically one can expect then that $y$ should be the closest trait to $x$ with a population large enough. However as the two traits are different, the individuals are not as compatible and the corresponding birth rate should decrease with the distance $|x-y|$.

The same phenomenon can be seen if one assumes that individuals of similar traits may cooperate. Selection-mutation models rightly focus on competition between individuals as the main interaction mechanism in order to observe selection (and hence evolution). Nevertheless cooperation usually exists as well; we assume here that it takes place at a smaller scale than the competitive interactions.

We are for example lead to add to the reproduction rate a cooperative effect of the kind 
\begin{equation}
\max\left(C_0,\quad \int_\R  \Phi\left(|x-y|/\eta_\ep\right)\;\frac{u_\ep(y)}{\eta_\ep}\,dy\right),   \label{correction0} 
\end{equation}
for a small parameter $\eta_\ep$ and a decreasing $\Phi$. Note that the maximum effect of cooperation is capped (obviously birth rates cannot blow up). 

Of course if one expects to see an effect on the small population then it is necessary that the reproduction rate without cooperation be negative, which means that the cooperation effect should be strong enough.

It is possible to simplify \eqref{correction0} (for numerical purposes for instance) while keeping the same structure and the same effects. When one combines it with the lower value of the reproduction rate, and assume for instance a linear $\Phi$, then it is essentially possible to replace \eqref{correction0} by
a regularization of 
\begin{equation}
-K\,d(x,\ \{u_\ep\geq \bar u_\ep\}.\label{correction1}
\end{equation} 
Those are the type of corrections that we consider here, with $K$ large enough.
\subsection{A brief overview of the various approaches for selection, mutation dynamics}
The stochastic approach is based on
individual-based models. As we mentioned before, they are related to 
evolutionary PDE models as those here or in~\cite{desvillettes-jabin-al-08,JR}
through a scaling of large 
population (see again \cite{CFM}). Using a simultaneous scaling of large
population and rare mutations, a stochastic limit 
process was obtained in~\cite{champagnat-06} in the case of a monotype
population (i.e.\ when the limit 
process can only be composed of a single Dirac mass), and in~\cite{CM}
when the limit population can be 
composed of finitely many Dirac masses. Other features can be added to those models, age-dependence for instance as in \cite{MT}. 

At the deterministic level, one approach consists in studying a simultaneous scaling of
mutation and selection, in order to 
obtain a limit dynamics where transitions from a single Dirac mass to
two Dirac masses could occur; the famous branching phenomenon that we have already mentioned. This is where the deterministic models presented here fit.
More is said about the contributions in that area in the next subsection.

Another approach consists in completely separating the two scalings. One then tries to directly characterize evolutionary dynamics as sums of
Dirac masses under biologically 
relevant parameter scalings, instead of obtaining it as some limit. This is a key point in \emph{adaptive
  dynamics}---\cite{HS,metz-nisbet-al-92,metz-geritz-al-96, 
dieckmann-law-96,CFBA}.

A classical way of justifying this form of the solution consists in studying 
the stationary behaviour of an evolutionary model involving a scaling
parameter for mutations, and then 
letting this parameter converge to 0. The stationary state has been
proved to be composed of one or several 
Dirac masses for various models (for 
deterministic PDE models, 
see~\cite{calsina-cuadrado-04,CCP,desvillettes-jabin-al-08,JR,GBV},
for Fokker-Planck PDEs corresponding to 
stochastic population genetics models, see~\cite{burger-bomze-96}, for
stochastic models, see~\cite{yu-07}, 
for game-theoretic models, see~\cite{cressman-hofbauer-05}). Closely
related to these works are the notions of 
ESS (evolutionarily stable strategies) and CSS (convergence stable
strategies)~\cite{metz-geritz-al-96,diekmann-04}, which allow one in
some cases to characterize stable 
stationary
states~\cite{calsina-cuadrado-04,
desvillettes-jabin-al-08,JR,cressman-hofbauer-05}.  

In this context the  
phenomenon of \emph{evolutionary branching}, which is especially important for us, simply corresponds to the direct
transition from a population composed of a single 
Dirac mass to a population composed of two Dirac masses, \cite{metz-geritz-al-96,geritz-metz-al-97,geritz-kisdi-al-98}. 
\subsection{An overview of the works around Eq. \eqref{u2}}
The scalings and the first formal results have been
obtained in~\cite{DJMP}. This was followed by several 
works on other models and on the corresponding
Hamilton-Jacobi 
PDE~\cite{CCP,PG}. 

The main difficulty as mentioned before is the identification of the weak limit of $u$ in terms of $\varphi$. There is usually one additional information which comes from uniform bounds on the total mass, it is that
\[
\max_\R \varphi(t,.)=0,\ \forall t.
\]
One solution is then to try to see the undetermined $I\star u$ as a Lagrange multiplier for this additional constraint. 

For example if $\varphi$ solves \eqref{phi0} then our constraint should imply that $r-I\star u$ is non positive on the set $\{\varphi=0\}$ and vanishes on at least one point.

Unfortunately, it is in general even formally not possible to identify $I\star u$ with that. For instance if $\varphi$ attains its maximum at one point then one has only constraint which is not enough.

There are however cases of competitive interactions where this suffices, Suppose for instance that $I=1$. Then $I\star u$ is just the total mass and formally one simply expects that
\[
\int_\R u(t,y)\,dy=\max_{x\;s.t.\;\varphi(t,x)=0} r(x).
\]
This is typical of so-called one resource interaction, meaning that the individuals interact only through one average quantity. In general that means considering interactions of the type
\begin{equation}
R\left(x,\ \int \eta(x)\,u(t,x)\,dx\right),\label{1resource}
\end{equation}
with an increasing $R$. 

For models of the type we consider here,
rigorous results (especially for the well posedness of the
Hamilton-Jacobi eq. at the limit) mainly only exist in this case with
just one resource, see \cite{BP} and \cite{BMP}, \cite{LMP} (one resource but
multidimensional traits).
 
However it is possible to extend the theory, see \cite{CJ}. Formally one expects the limit $u$ to satisfy the following conditions
\begin{equation}
\begin{split}
& (i)\ \mbox{supp}\,u(t,.)\subset \{\varphi(t,.)=0\},\\
& (ii)\ r-I\star u\leq 0,\quad \mbox{on}\ \{\varphi(t,.)=0\},\\
& (iii)\ r-I\star u= 0,\quad \mbox{on}\ \mbox{supp}\,u(t,.).\\
\end{split}\label{ESS}\end{equation}
This corresponds to the definition of an evolutionarily stable strategy. If one assumes a strong competition, {\em i.e.} that the operator $I\star $ is positive, then there exists a unique measure satisfying \eqref{ESS} (see \cite{JR}). 

This is the approach followed here as well. With respect to the case with one resource, there are nevertheless additional difficulties. It is much harder to control the time oscillations of the reaction term, which is usually increasing and hence $BV$ in time for one resource.
\section{The model studied and the results}%
According to the previous considerations, we study the following equation
\begin{equation}
\partial_t u_\eps(t,x)=\frac{1}{\eps}\left(r-I\star_x u_\ep(t,.)-D_\eps(u_\eps)\right)\,u_\eps(t,x)
+M_\eps(u_\eps)(t,x), \label{ueps} 
\end{equation}
where we recall that the mutation kernel $M_\eps$ reads
\begin{equation}
M_\eps(f)(x)=\frac{1}{\eps}\int_\R K(z)\,(f(x+\eps z)-f(x))\,dz,\label{Meps}
\end{equation}
for a $K\in C^\infty_c(\R)$ such that $\int_\RR zK(z)\,dz=0$. 
Following \cite{DJMP}, one defines $\varphi_\eps$ as 
\begin{equation}
u_\eps=e^{\varphi_\eps/\eps},\quad \mbox{or}\ \varphi_\eps=\eps\,\log u_\eps,
\end{equation}
and obtains the equation
\begin{equation}
\partial_t \varphi_\eps=r-I\star u_\ep+H_\eps(\varphi_\eps),\label{phieps}
\end{equation}
with
\begin{equation}
H_\eps(f)=\int_\R
K(z)\,\left(e^{(f(x+\eps\,z)-f(x))/\eps}-1\right)\,dz.\label{Heps} 
\end{equation}
First of all we need to make sure that only the traits in a compact interval are important (to avoid traveling waves effects for instance). This can be simply ensured by asking all traits out of an interval to have a negative reproduction rate
\begin{equation}
\exists R>0,\quad \exists r_0>0, \quad \forall |x|>R,\quad r(x)\leq -r_0.
\label{boundeta}
\end{equation}
Of course the whole population should not vanish immediately, and it is necessary that a non negligible part be concentrated on $[-R,\ R]$. The total population should also be bounded which leads to
\begin{equation}\begin{split}
&\sup_\ep \int_{\R} u_\ep(t=0,.)\,dx<\infty,\\
& \exists R_0\,\ \min_{|x|<R_0} r(x)>0,\qquad \inf_\ep \int_{-R_0}^{R_0} u_\ep(t=0,x)\,dx>0.\end{split} \label{boundpop}
\end{equation} 
We assume that the individuals interact through a strong competition
\begin{equation}
\forall f\in M^1(\R)\setminus\{0\},\quad \int_{\R^2} I(x-y)\,df(x)\,df(y)>0,\qquad \mbox{or}\ \hat I>0.\label{strongcompet}
\end{equation}
This allows us to define uniquely the ESS as per
\begin{prop} Assume \eqref{strongcompet}, \eqref{boundeta} and that $r,\;I\in C(\R)$. For any closed $\Omega\subset \R$, there exists a unique
  finite nonnegative measure $\mu(\Omega)$ satisfying\\ 
  i) $\mbox{supp}\,\mu\subset \Omega$\\
  ii) 
$r-I\star \mu\leq 0\ in\ \Omega$,
$\quad r-I\star\mu = 0\ on\ \mbox{supp}\,\mu$.\label{metastable}
\end{prop}
Note however that sometimes, one may have uniqueness of the
environmental variables whereas the population measure is not unique
(and \eqref{strongcompet} is violated), most of our method would
remain valid in such a case. This is the situation of a single
resource, where almost nothing is required, see~\cite{BP}.
 Nevertheless in
more general situations, the conditions for which this kind of property
holds are not currently identified.

We need an additional assumption to make sure that the ESS can only be concentrated on a set of measure $0$
\begin{multline}
\exists S\in C(\R_+)\ \mbox{with}\ S(0)=0\ s.t.\ \forall \mu\in M^1(\R_+)\\ |\{x,\ |r(x)-I\star_x \mu(x)|\leq \nu\}|\leq S(\nu).\label{rootsnumber} 
\end{multline}
This condition is in part technical and is required to avoid some time oscillations. However it also corresponds to the natural biological idea that only a few traits can be present at a given time.

It is probably hard to check~\eqref{rootsnumber} in specific
models, but it is at least satisfied in 
large classes of parameters. One easy example is if the derivatives
$r^{(k)}$ and $I^{(k)}$ are positive (or negative) for some $k$. It has in fact been proved, see \cite{Gy}, that generically in $r$ and $I$ the ESS is discrete (a finite sum of Dirac masses) and hence \eqref{rootsnumber} should be satisfied. 

As for $D_\ep$, we assume that there exists a critical scale $\fe$ s.t.
\begin{equation}\begin{split}
&\sup_\ep \sup_{t,x} (|D_\ep|+|\partial_x D_\ep(t,x)|)<\infty,\quad \sup_\ep \sup_{t,x} \fe\,|\partial_{xx} D_\ep(t,x)|<\infty,\\
& D_\ep(t,x)=0\quad\mbox{if}\quad \pe(t,x)\geq -\fe\\
& |\partial_x D_\ep|\geq 2\,|\partial_x r|\quad \mbox{if}\quad d(x,\ \{y,\ \pe(t,y)\geq -\fe\})\geq 2\,\fe.\\
\end{split}\label{assdep}
\end{equation}
Those assumptions are compatible with the type of corrections like \eqref{correction0} and \eqref{correction1}. They introduce a new scale in the problem $-\fe$ which corresponds to a critical population density of $\exp(-\fe/\ep)$. In line with our previous consideration, we would like to take this polynomial in $\ep$ which means $-\fe\sim \ep\log\ep$ or more precisely for a constant $C$ uniform in $\ep$ 
\begin{equation}
\ep\,\log \frac{1}{\ep}\leq \fe\leq C\,\ep\,\log\frac{1}{\ep}.\label{sizefe} 
\end{equation}
There is no need to be more precise on $D_\ep$ to study the properties of Eq. \eqref{phieps}. However if one wants to identify the limit, it is necessary to specify what $D_\ep$ should look like at the limit. Therefore we make the additional assumption, for any fixed $f\in C(\R)$, non positive
\begin{equation}
D_\ep=D_\ep(f)\longrightarrow \min (K\,d(x,\ \{f=0\}),\ D_0)\qquad
\mbox{as}\quad \ep\rightarrow 0.\label{limitdep}
\end{equation}
This corresponds to \eqref{correction1} but other shapes would be possible and would essentially work the same.

Formally, we can hence expect that as $\eps\rightarrow 0$, Eq. \eqref{phieps} will lead to
\begin{equation}
  \label{eq:closed-HJ}
  \partial_t\varphi=r-I\star_x\mu(\{\varphi(t,.)=0\})-\min (K\,d(x,\ \{f=0\}),\ D_0)+H(\partial_x\varphi). 
\end{equation}
where the Hamiltonian $H_\eps$ became
\begin{equation}
H(p)=\int_\RR K(z)\,\left(e^{p\,z}-1\right)\,dz.\label{defH}
\end{equation}
This is indeed what one can prove
\begin{theorem} Assume $K,\;r\in L^\infty\cap C^{2}_c(\RR)$, $\int_\RR
  zK(z)\,dz=0$, \eqref{boundeta} on $r$, \eqref{boundpop} on the bounds for the initial population, \eqref{strongcompet} on the interaction $I$
  \eqref{rootsnumber}, \eqref{assdep} on $D_\ep$ with \eqref{sizefe}, and \eqref{limitdep} on the limit of $D_\ep$, that the initial data
  $u_\eps(t=0)>0$ or $\pe(t=0)$ are 
  $C^2$, with
  \begin{gather}
    \inf_\eps\inf_{x\in\RR}\partial_{xx}\varphi_\ep(t=0,x)>-\infty,
  \end{gather}
  and that $\varphi_\varepsilon(t=0,\cdot)$ converges to a function
  $\varphi^0$ for the norm 
  $\|\cdot\|_{L^{\infty}(\R)}$.

  Then up to the extraction of a subsequence in $\eps$, $\pe$ converges to some
  continuous $\varphi$ uniformly on any compact subset of 
  $[0,T]\times\R$ and $\varphi$ is a solution to 
  \eqref{eq:closed-HJ} almost everywhere in $t,x$ with initial condition
  $\varphi(t=0,\cdot)=\varphi^0$. In particular the function $I\star u_\ep$
  converges to $I\star u$ in $L^p([0,T],\ C(\R))$ for any 
  $p<\infty$, where $u(t,.)=\mu(\{\varphi(t,.)=0\}$ is defined from $\varphi$ by Prop. \ref{metastable} and is continuous in time. 
  \label{theolimit}
\end{theorem}

From a practical point of view, computing the solution $u_\eps$ of
Eq. \eqref{ueps} is often too costly for small $\eps$. This result
allows 
to approximate
the population density $u_\eps$ for 
small $\varepsilon$ by the simpler $\mu(\{\varphi(t,\cdot)=0\})$,
where $\varphi$ may be obtained by a 
discretization of~(\ref{eq:closed-HJ}), in the fashion of those done
in~\cite{DJMP}. Rigorous numerical 
 analysis of this kind of Hamilton-Jacobi equations is however 
still very preliminary. 

\bigskip

Section \ref{secprop} gives a short sketch of the proof of Prop. \ref{metastable} and can be safely ignored if one is familiar with \cite{JR}, \cite{Ra}.

The proof of Theorem \ref{theolimit} is given in the next section and being quite technical is divided into several lemmas. Lemma \ref{apriori} just corresponds to the classical a priori estimates on equations like \eqref{phieps}. Lemma \ref{conttime} essentially follows the step of \cite{CJ} and can also be skipped if the reader is familiar with that work. 

In the proofs below, $C$ denotes a numerical constant which may change
from line to line but which only depends on $T$, norms of the initial data or the coefficients, but which is always uniform on $\ep$.  

We define as usual the distance between a point and a set
\[
d(x,\Omega)=\inf_{y\in \Omega} |x=y|.
\]
We also define the semi distance $\delta(O_1,O_2)$ between two sets $O_1$ and $O_2$ as usual by
\[
\delta(O_1,O_2)=\sup_{x\in O_1}\inf_{y\in O_2} |x-y|.
\] 
A small $\delta(O_1,O_2)$ indicates that $O_1$ is almost included in $\bar O_2$ and in particular $\delta(O_1,O_2)=0$ iff $O_1\subset \bar O_2$.

We denote by $M^1(\omega)$ the set of signed Radon measures on the
subset $\omega$ of $\R$ equipped with the 
total variation norm.
%
%
\section{Sketch of the proof of Prop. \ref{metastable}}\label{secprop}
A complete proof of this proposition can be found in \cite{JR}, and \cite{Ra}. We only show uniqueness and give a rough skecth of the existence part. For a given compact set $\Omega$, assume that we can find two measures $\mu_1$ and $\mu_2$ with support in $\Omega$ and such that
\[
r-I\star\mu_i=0\quad\mbox{on}\ \mbox{supp}\,\mu_i,\qquad r-I\star \mu_i\leq 0\quad \mbox{on}\ \Omega.
\]
Compute
\[
0\geq \int (r-I\star \mu_1)\,d\mu_2=\int I\star (\mu_2-\mu_1)\,d\mu_2.
\]
Summing with the symmetric term, one obtains
\[
\int I\star(\mu_2-\mu_1)\,(d\mu_2-d\mu_1)\leq 0.
\]
However since $\hat I>0$ then the corresponding quadratic form is positive and we can conclude that $\mu_1=\mu_2$.

\medskip

For the existence, one considers the equation
\[
\partial_t \mu_\ep=\frac{1}{\ep} (r-I\star \mu_\ep)\,\mu_\ep+\ep\,M(u_\ep),\quad \mu_\ep(t=0)=\mu^0,
\]
for a well chosen initial data $\mu^0$. If $\Omega$ is fully discrete, $\Omega=\{x_1,...,x_n\}$ then one may choose $\mu^0=\sum_i \delta_{x_i}$. If $\Omega$ is an interval then one simply takes $\mu_\ep=1$. The general case is trickier and consists in taking an initial measure $\mu^0$ s.t. for any $\delta>0$ and $x_0\in \Omega$
\[
\int_{\Omega\cap [x_0-\delta,\ x_0+\delta]} d\mu^0>0.
\]
It is then easy to obtain lower and upper bounds for the total mass.

Note that though there is also a mutation term here, its scaling is completely different from \eqref{ueps}. The scaling does not involve small mutations and moreover the mutations are of order $\ep$ (instead of $1$) thus vanishing at the limit.

Passing to the weak limit, $\mu_\ep\rightarrow \mu$, gives a measure that has the desired properties (but is by no means easy to show).
%
\section{Proof of Theorem \ref{theolimit}}

\subsection{A Priori estimates}
We start by stating and proving the obvious a priori estimates that one can obtain for the problem. Those essentially follow the lines of previous works.

We show the following estimates on the solution to \eqref{phieps}
\begin{lemma} Let $\varphi_\eps$ be a solution to \eqref{phieps} with
  the assumptions of Theorem \ref{theolimit}. Then for any $T>0$ 
  \begin{eqnarray}
  &\|\partial_t\varphi_\eps\|_{L^\infty([0,T]\times\R)}+\|\partial_x
  \varphi_\eps\|_{L^\infty([0,T], 
      L^\infty(\R))}\leq C_T,\\
 &   \forall\,t\leq T,\,x\in \R,\quad \partial_{xx}\pe(t,x) \geq -\frac{C}{\fe},\quad
    H_\eps(\pe)\geq -C \frac{\eps}{\fe}, \label{infHeps}\\
  &  \forall\,t\leq T,\ \frac{1}{C}\leq \int_\R u_\eps(t,x)\,dx\leq C,\quad
    \pe(t,x)\leq \eps\,\log 1/\eps+C\,\ep, 
  \end{eqnarray}
 where $C$ only depends on the time $T$, $\int_\mathbb{R}
  e^{\varphi_\varepsilon^0(x)/\varepsilon}\,dx$, 
  $\|\partial_x\varphi_\varepsilon^0\|_{L^\infty(\R)}$ and the infimum of
  $\partial_{xx}\varphi_\varepsilon^0(x)$. Finally $\pe$ has level sets, uniformly bounded in $\ep$. \label{apriori}
\end{lemma}

\noindent{\bf Proof.}  Almost all proofs here are taken directly from \cite{CJ} and reproduced for the sake of completeness. Some resulting bounds are different and much worse than in this former article, more precisely the lower bounds on $\partial_{xx} \pe$ and $H_\ep$. Nevertheless even for those, the proofs, {\em i.e.} the way to obtain the bounds, are very close. As such we may skip some technical details. 

\medskip

{\em Step 0: Upper Bound on the total mass.} 
First
notice that because of \eqref{boundeta}, there exists $R>0$ s.t.
\[
\forall |x|>R,\quad r(x)-I\star u_\ep(t,x)\leq -r_0.
\]
Let $\psi$ be a smooth test function with support in $|x|>R$, taking
values in $[0,1]$ and equal to $1$ on 
$|x|>R+1$. Using the previous bound, we compute
\[\begin{split}
\frac{d}{dt}\int_\R \psi(x)\,u_\eps(t,x)\,dx\leq
&-\frac{r_0}{\eps}\int_\R
  \psi(x)\,u_\eps(t,x)\,dx\\
&+\frac{1}{\eps}\int_{\R^2} K(z) (\psi(x-\eps
  z)-\psi(x))\,u_\eps(t,x)\,dz\,dx\\
&\leq -\frac{r_0}{\eps}\int_\R
  \psi(x)\,u_\eps(t,x)\,dx+C\int_\R u_\eps(t,x)\,dx.  
\end{split}\] 

On the other hand, on the bounded domain $[-R-1,\ R+1]$ as $I(x)>0$ for all $x$,  one has for some constant $C$
\[
\forall |x|<R+1,\quad I\star u_\ep(t,x)\geq C\int_\R (1-\psi)\,u_\ep\,dx.
\]
Therefore with the same kind of estimate
\[\begin{split}
\frac{d}{dt}\int_\R (1- \psi(x))\,u_\eps(t,x)\,dx\leq
&C\int_\R u_\eps(t,x)\,dx\\
+\frac{1}{\eps} \Bigg(\sup r-&C\,\int_\R (1-\psi)\,u_\eps\,dx\Bigg)\,\int_\R (1-\psi(x))\,u_\eps(t,x)\,dx.\\
\end{split}\]
Summing the two
\[\begin{split}
\frac{d}{dt}\int_\R u_\eps(t,x)\,dx\leq&\frac{1}{\eps} \left(\sup r-C\,\int_\R (1-\psi)\,u_\eps\,dx\right)\,\int_\R (1-
\psi)\,u_\eps\,dx \\
&-\frac{r_0}{\eps}\int_\R
  \psi(x)\,u_\eps(t,x)\,dx+C\int_\R u_\eps(t,x)\,dx.
\end{split}\] Since the sum of the first two terms of the r.h.s.\ is
negative if $\int u_\varepsilon$ is 
larger than a constant independent of $\varepsilon$, this shows that
$\int u_\eps(t,x)\,dx$ remains uniformly 
bounded on any finite time interval.


\medskip
 
{\em Step 1: Bound on $\partial_x \pe$.} This is a classical bound for
solutions to Hamilton-Jacobi equations and  we have to
check that it remains true uniformly at the $\eps$ level. We follow exactly \cite{CJ} for instance. Compute
\begin{equation}\begin{split}
&\partial_t \partial_x \pe=\partial_x r-\partial_x I\star u_\ep-\partial_x D_\ep\\
&\ +\int
K(z)\,
e^{\frac{\pe(t,x+\eps z)-\pe(t,x)}{\eps}}\,\frac{\partial_x\pe(t,x+\eps
  z)-\partial_x\pe(t,x)}{\eps}\,dz. \label{eq:dtx} 
\end{split}\end{equation}
By our assumptions and the upper bound on the total mass
\[
|\partial_x r-\partial_x I\star u_\ep-\partial_x D_\ep|\leq C.
\]
First note that this shows that $\|\partial_x\pe(t,.)\|_{L^\infty}$ remains finite over a (possibly very short) time interval $[0,\ t_\ep]$.

Now we use the classical maximum principle. Fix $t\in[0,T]$ such that
$C_{\varepsilon,t}:=\|\partial_x\pe(t,\cdot)\|_{L^\infty(\R)}<\infty$. For
any $x\in\R$ such that 
$\partial_x\varphi_\varepsilon(t,x)>\sup_y
\partial_x\varphi_\varepsilon(t,y)-\alpha$, where the 
constant $\alpha>0$ will be specified later, we have
$$
\partial_t\partial_x\varphi_\varepsilon(t,x)\leq C+\int_\R K(z) 
e^{|z|\,C_{t,\varepsilon}}\frac{\alpha}{\varepsilon}\,dz\leq
C\Big(1+\frac{\alpha}{\varepsilon}e^{C\,C_{t,\varepsilon}}\Big).
$$
Therefore, choosing $\alpha=\varepsilon e^{-C\,C_{t,\varepsilon}}$, we obtain
\[
\frac{d}{dt}\, \sup_x\partial_x\pe(t,x)\leq C,
\]
for a constant $C$ independent of $t<t_\varepsilon$ and of
$\varepsilon$. Using a similar 
argument for the minimum, we deduce that 
$t_\varepsilon>T$ and that $\partial_x\varphi_\varepsilon$ is bounded
on $[0,T]\times\R$ by a constant 
depending only on $T$ and $\|\partial_x\varphi_\varepsilon^0\|_{L^\infty(\R)}$.

\medskip

{\em Step 2: First bound on $H_\eps(\pe)$ and bounds on 
$\partial_t\pe$ and $\pe$.} First remark that
\[\begin{split}
-\int_\R K(z)\,dz\leq H_\eps(\pe(t))(x)&=\int_\R
K(z)\,
e^{\frac{\pe(t,x+\eps z)-\pe(t,x)}{\eps}}\,dz-\int_\R K(z)\,dz\\
&
\leq \int
K(z)\,e^{|z|\,\|\partial_x \pe\|_{L^\infty([0,T],\R)}}\,dz\leq C.
\end{split}\]
This is not optimal, more precisely the lower bound is atrocious, but will suffice for the moment. 

Directly from Eq. \eqref{phieps},
\[
|\partial_t \pe|\leq \sup |r|+\sup |I|\,\int_\R u_\ep\,dx+|D_\ep|+C\leq C,
\] 
hence ending the proof of the full Lipschitz bound on $\pe$.

To get the upper bound on $\pe$, we use this Lipschitz bound to get
\[
\pe(t,y)\geq \pe(t,x)-C\,|y-x|,
\]
so
\[
\int_\R u_\eps(t,y)\,dy\geq \int_\R e^{\pe(t,x)/\eps}\,e^{-C\,
|y-x|/\eps}
\,dy\geq 2 C^{-1}\,\eps\,e^{\pe(t,x)/\eps}.
\]
Hence the bound on the total mass yields that $\pe\leq \,\eps\,\log
1/\eps+C\,\ep$. 

\medskip

{\em Step 3: Lower bounds on $\partial_{xx} \pe$ and $H_\ep$.} Here we start paying for the introduction of an additional mortality term as the second derivative of $D_\ep$ is not bounded uniformly in $\ep$. As before we use a maximum principle, from \eqref{phieps} and \eqref{assdep}
\[\begin{split}
\partial_t \partial_{xx}\pe&\geq -\frac{C}{\fe}+\int_\R
K(z)\,
e^{\frac{\pe(t,x+\eps z)-\pe(t,x)}{\eps}}\,\frac{\partial_{xx}\pe(t,x+\eps
  z)-\partial_{xx}\pe(t,x)}{\eps}\,dz\\
& +\int_\R
K(z)\,
e^{\frac{\pe(t,x+\eps z)-\pe(t,x)}{\eps}}\,\frac{(\partial_{x}\pe(t,x+\eps
  z)-\partial_{x}\pe(t,x))^2}{\eps}\,dz. 
\end{split}\]
The last term is of course non negative and so with the same argument
as before, we get
\[
\frac{d}{dt} \inf_x \partial_{xx}\pe(t,x)\geq -\frac{C}{\fe}.
\]
This proves the uniform
lower bound on $\partial_{xx} \pe$. Let us turn to the sharp lower bound on 
$H_\eps(\pe)$. Let us write
\[
H_\eps(\pe)\geq \int_\R K(z)\,\exp\left(\int_0^1
z\,\partial_x\pe(t,x+\theta z\,\eps)\,d\theta\right)\,dz-\int_\R K(z)\,dz.
\]
By differentiating once more
\[\begin{split}
\int_0^1
z\,\partial_x\pe(t,x+\theta
z\,\eps)\,d\theta&\geq z\,\partial_x\pe(t,x)\\
&\qquad+\int_0^1z \int_0^1
\theta\,z\,\eps\,
\partial_{xx}\pe(t,x+\theta'\theta z\eps)\,d\theta'd\theta\\
&\geq z\,\partial_x\pe(t,x)-C\,\frac{\eps}{\fe}\,z^2.
\end{split}\]
Eventually
\begin{align*}
  H_\eps(\pe) & \geq \int_\R K(z) \exp(
  z\,\partial_x\pe(t,x)-C\,\eps\,z^2/\fe)\,dz-\int_\R K(z)\, dz \\ & \geq
  H(\partial_x\pe(t,x))-C\,\frac{\eps}{\fe},
\end{align*}
where $H$ is defined as in \eqref{defH} and since $K$ is compactly
supported. Because we assumed that $\int_\R zK(z)\,dz=0$, we have
$H(p)\geq 0$ for any $p$, which gives the final bound.

\medskip

{\em Step 4: $\pe$ has uniform compact level sets.}\\
Observe that $\varphi_\varepsilon(t=0,x)\rightarrow-\infty$ when
$x\rightarrow\pm\infty$ since $\int_\R 
u_\varepsilon(t=0,x)\,dx<\infty$ and $\partial_x\varphi_\ep(t=0)$ is
bounded. Because of the uniform convergence 
of $\varphi_\varepsilon(t=0)$ to $\varphi^0$ on $\R$, one moreover deduces that
the convergence $\varphi_\varepsilon(t=0,x)\rightarrow-\infty$ is uniform in $\ep$.

Since $\partial_x\varphi_\ep \in L^\infty([0,T],\R)$,
this remains true at any time $t\in [0,\ T]$, uniformly in $\ep$.

Therefore, the set
$$
\Omega_\ep:=\{(t,x)\in[0,T]\times\R:\varphi_\ep(t,x)\geq -l\}
$$
is bounded for any $l$, uniformly in $\ep$.

\medskip

{\em Step 5: Lower bound on the total mass and $\pe$.}\\
By the previous steps we know that 
\[
M_\ep(u_\ep)\geq \frac{u_\eps}{\eps} H_\eps\geq -C\,\frac{u_\ep}{\fe}.
\]

Note that $D_\ep=0$ whenever $u_\ep\geq e^{-\fe/\ep}$. 
Therefore by \eqref{boundpop}, choose $R_0$ s.t. $\min_{|x|\leq K} r>0$ and  integrate \eqref{ueps}
\[
\frac{d}{dt} \int_{|x|\leq K} u_\ep\,dx\geq \frac{1}{\ep}\left(\min _{|x|\leq K} r-C\,\frac{\ep}{\fe}-C\,\int_\R u_\ep\,dx\right)\, \int_{|x|\leq K} u_\ep\,dx.
\]
This implies that by \eqref{boundpop}
\[
\log \int_{|x|\leq K} u_\ep(t,x)\,dx\geq -C+ \frac{1}{\ep}\int_0^t \left(\min _{|x|\leq K} r-C\,\frac{\ep}{\fe}-C\,\int_\R u_\ep(s,x)\,dx\right)\,ds,
\]
or obviously
\[
\log \int_{\R} u_\ep(t,x)\,dx\geq -C+ \frac{1}{\ep}\int_0^t \left(\min _{|x|\leq K} r-C\,\frac{\ep}{\fe}-C\,\int_\R u_\ep(s,x)\,dx\right)\,ds.
\]
This allows us to conclude that the total mass remains bounded from below uniformly in $\ep$. 

In particular by step 4, this lower bound means that $\max \pe\geq -\fe$.
%
%
\subsection{Passing to the limit in the equation: First steps}
%
By the uniform bounds provided by Lemma \ref{apriori}, we can extract a subsequence in $\eps$ (still denoted with
$\varepsilon$), and find a function $\varphi$ 
on $[0,T]\times\R$ such that $\partial_t\varphi\in
L^\infty([0,T]\times\R)$, $\partial_x \varphi\in 
L^\infty([0,T]\times\R)$, with $\max \varphi=0$
satisfying by the Arz{\'e}la-Ascoli theorem
\begin{equation}\begin{aligned}
\pe & \longrightarrow \varphi 
\quad\mbox{uniformly in\ }C(K)\mbox{\ for any
      compact\ }K\mbox{\ of\ }[0,T]\times\R,\\ 
\end{aligned}\label{limitphieps}
\end{equation}
Note that by the upper and lower bounds on $\pe$, one has $\max \varphi=0$.

Since $u_\ep$ is uniformly bounded in $L^1$, it converges (still after an extraction) in the weak-* topology of measures to some $u\in L^\infty([0,\ T],\ M^1(\R))$. This in turn implies a weak convergence of $I\star u_\ep$. Note that it is indeed only a weak convergence in spite of the convolution because it regularizes only in $x$ and time oscillations are still possible.

Using the notion of viscosity solution (see \cite{BP} for instance) and the uniform bound on $\partial_x \pe$, we may obtain the convergence of $H_\ep(\pe)$ to $H(\partial_x \varphi)$.

Finally as $D_\ep$ is uniformly bounded and can hence be assumed to converge weakly to some $D\in L^\infty$, we obtain
\begin{equation}
\partial_t \varphi=r-I\star u(t,.)-D+H(\partial_x\varphi).\label{limeq}
\end{equation}
Unfortunately $D$ and $u$ are still unidentified. This is now our aim but it will require a much more precise understanding and control of the time oscillations of the set where the population is concentrated.
%
%
%
\subsection{Continuity in time of the set $\{x\,|\; \varphi=0 \}$}
At the limit, it is possible to show that the points where the
population is concentrated move at most at a finite speed given by
\[
V=2\,\sup_\eps\,\sup_{|\xi|\leq 2\|\partial_x\pe\|_{L^\infty}}\,\frac{1}{|\xi|}
\int e^{z\,\xi}\,K(z)\,dz-1.
\]  
We comment on that again at the very end of the proof. However we do not have yet enough tools to prove such a strong statement.
For the time being we will be satisfied in proving a continuity result. The first step is to do it at the limit
\begin{lemma}
Assume that for some point $x_0$, some $t>t_0$ and $\delta>0$ 
\[
d(x_0,\,\{\varphi(s,.)=0\})\geq \delta,\qquad \forall\ s\in[t_0,\ t] ,
\]
then $\varphi(t,x_0)<0$.\label{continuitylimit}
\end{lemma}

Before turning to the next result, let us point out that Lemma \ref{continuitylimit} actually implies the following
\begin{equation}\begin{split}
&\exists \tau\in C(\R_+)\ \mbox{with}\ \tau(0)=0,\ s.t.\ \forall s\geq t,\\
&\forall x\in \{\varphi(s,.)=0\},\quad \exists y\in \{\varphi(t,.)=0\}\ \mbox{with}\ |y-x|\leq \tau(s-t).
\end{split}\label{deftau}
\end{equation} 

One of the main idea in the following proofs is to follow the characteristics
corresponding to the Hamiltonian $H$. However Eq. \eqref{phieps}
involves the modified Hamiltonian $H_\ep$ which does not have
characteristics per se (it is non local for example). Some
modifications are hence needed and follow the usual ideas 
for parabolic problems, which is why we use the following lemma

\begin{lemma}
Introduce the intermediary scale $\sqrt{\eps\,|\fe|}=\eps/\bar\eps=\bar\eps\,\fe$ and consider any interval $[a,\ b]$. Define $a(t)=a+(t-t_0)\,V/2$, $b(t)=b-(t-t_0)\,V/2$. We denote
\[
m_\eps(t)=\max_{[a(t),\ b(t)]} \pe(t,.).
\] 
Then
\[\begin{split}
&m_\eps(t)-m_\eps(t_0)\geq \int_{t_0}^t\Big(\min_{[a(t),\ b(t)]} (r(x)-I\star_x u_\eps(s,x)-D_\eps)-\frac{C}{\eps}{\fe}\Big)\,ds \\
&m_\eps(t)\leq \max_{[a-\bar\ep\fe,\ b+\bar\ep\fe]} \pe(t_0,.)\\
&\qquad\quad+\int_{t_0}^t\max_{[a(s)-\bar\eps\,\fe,\ b(s)+\bar\eps\,\fe]} (r(x)-I\star_x u_\eps(t,x)-D_\eps+2\,\bar\ep)\,ds. 
\end{split}\]\label{minmaxlemma}
\end{lemma}
\noindent{\bf Proof of Lemma \ref{minmaxlemma}.} The part
\[
\min_{[a(t),\ b(t)]} (r(x)-I\star_x u_\eps(s,x)-D_\ep)-\frac{C}{\eps}{\fe} \leq\frac{d}{dt} m_\eps(s)
\]
is direct once one recalls that for any $x\in [a(t),\ b(t)]$ then
\[\begin{split}
\partial_t \pe(s,x)&=r(x)-I\star_x u_\eps(s,x)-D_\eps+H_\eps(\pe)(s,x)\\
&\geq \min_{[a(t),\ b(t)]} (r(x)-I\star_x u_\eps(s,x)-D_\eps)-\frac{C}{\eps}{\fe},
\end{split}\]
by the lower bound \eqref{infHeps}.

Let us turn now to the upper bound which is trickier and involves the velocity $V$.

Consider $\chi_\eps$ a regularization of
$\frac{3}{2}\|\partial_x\pe\|_{L^\infty}\,()_+$, with $()_+$ the positive part,
at the intermediary scale $\bar\eps\,\fe$
, {\em i.e.}
\begin{equation}\begin{split}
&\chi_\eps(x)=0,\quad\mbox{if}\ x\leq 0,\\
&\chi'_\eps(x)=\frac{3}{2}\|\partial_x\pe\|_{L^\infty},\quad\mbox{if}\ x\geq \bar\eps\fe,
\end{split}\label{basicchi}
\end{equation}
together with the general bounds
\begin{equation}\begin{split}
& 0\leq \chi_\eps\leq
\frac{3}{2}\|\partial_x\pe\|_{L^\infty}\,(x)_++\bar\eps\,\fe, 
\quad \chi_\eps\geq
\frac{3}{2}\|\partial_x\pe\|_{L^\infty}\,(x)_+-\bar\eps\,\fe,\\ 
&0\leq \chi_\eps'(x)\leq 2\,\|\partial_x\pe\|_{L^\infty}\,
\ind_{x\geq 0}+\bar\eps\,\fe,\quad 
|\chi_\eps''(x)|\leq \frac{C}{\eps}\,\bar\eps.
\end{split}\label{propchi}
\end{equation}
With that one defines
\[
\psi_\ep(t,x)=-\chi_\eps(a(t)-x)-\chi_\eps(x-b(t)).
\]
Note that $\psi_\eps$
satisfies the following inequation
\[
\partial_t \psi_\eps \leq
-\frac{V}{2}\,|\partial_x\psi_\eps|.
\]
This has for consequence that $\forall z\in \mbox{supp}\,K\subset [-1,\ 1]$
\begin{equation}\begin{split}
\int_{\R}\exp&\left(
\frac{|\psi_\eps(t,x+\eps\,z)-\psi_\eps(t,x)|}{2\,\eps}\right)
\,K(z)\,dz-1\\
&\qquad \leq
\int_{\R}e^{z\,\partial_x\psi_\eps(t,x)}\,K(z)\,dz-1
+\bar\eps\leq \frac{V}{2}\,|\partial_x\psi_\eps|+\bar\eps.\label{ineqpsi}  
\end{split}\end{equation}

As $\pe$ and $\psi_\ep$ are
smooth functions, it is easy to study
\[
\tilde m_\eps(t)=\sup (\pe(t,.)+\psi_\ep(t,.))=\max (\pe(t,.)+\psi_\ep(t,.)).
\]
In general this maximum is attained at one (or several) point $x\in
\omega_m(t)$. First note that at such a point, one has
\[
\partial_x (\pe(t,x)+\psi_\ep(t,x))=0.
\]
By \eqref{basicchi} this implies that
\begin{equation}
\omega_m(t)\subset [a(t)-\bar\eps\fe,\ b(t)+\bar\eps\fe].\label{inclusionomegam}
\end{equation}
Now in general
\[
\frac{d}{dt} \tilde m_\ep(t)\leq \sup_{x\in \omega_m(t)} \partial_t (\pe(t,x)+\psi_\ep(t,x)).
\]
Note that by the definition of the maximum, $\forall z\in \mbox{supp}\,K$
\[
\pe(t,x+\eps\,z)\leq \pe(t,x)+\psi_\eps(t,x)-\psi_\eps(t,x+\eps\,z).
\]
Now for $x\in \omega_m(r)$
\[\begin{split}
H_\eps(\pe)(t,x)&= \int_{\R}
\exp\left(\frac{\pe(t,x+\eps\,z)-\pe(t,x)}{\eps}\right)
\,K(z)\,dz-1 \\
&\leq \int_{\R}
\exp\left(\frac{-\psi_\eps(t,x+\eps\,z)+\psi_\eps(t,x)}{\eps}\right)
\,K(z)\,dz-1\\
&\leq \frac{V}{2}\,|\partial_x\psi_\ep|+\bar\eps\leq 
-\partial_t \psi_\ep+\bar\eps. 
\end{split}\]
Therefore
\[\begin{split}
\frac{d}{dt} \tilde m_\eps(t)&\leq \sup_{x\in\omega_m(t)} (r(x)-I\star_xu_\ep(t,x)-D_\ep)+\bar\eps\\
&\leq \sup_{x\in [a(t)-C\,\bar\eps\fe,\ b(t)+C\,\bar\eps\fe]} (r(x)-I\star_x u_\eps(t,x)-D_\ep) + \bar\eps,\\
&\leq \sup_{x\in [a(t),\ b(t)]} (r(x)-I\star_x u_\eps(t,x)-D_\ep) + 2\,\bar\eps
\end{split}\]
by \eqref{inclusionomegam}.

Note that by \eqref{basicchi} and \eqref{inclusionomegam}
\[
m_\eps(t)\leq \tilde m_\eps(t)\leq \max_{[a(t)-\bar\ep\fe,\ b(t)+\bar\ep\fe]} m_\eps(t),
\]
which allows to conclude.
 \cqfd

Let us now turn to  Lemma \ref{continuitylimit}

\noindent{\bf Proof of Lemma \ref{continuitylimit}.} 

{\em Step 1: Semi-continuity of $\{\varphi=0\}$}\\
Let us start with the following crucial observation
\begin{equation}\begin{split}
 & \forall t,\ \exists \tau_{t}\in
  C(\R_+)\ \mbox{with}\ \tau_{t}(0)=0,\ \mbox{s.t.}\ \forall s\geq
  t,\ \\ 
&\forall 
  x\in \{\varphi(s,.)=0\},\;\exists y\in
  \{\varphi(t,.)=0\}\ \mbox{with}\ |y-x|\leq \tau_t(s-t). 
\end{split}\end{equation}
This is a sort of semi-continuity for $\{\varphi=0\}$. It is proved
very simply by contradiction. If it were 
not true, then
\[\begin{split}
&\exists t,\,\exists \tau_0>0,\,
\exists s_n\rightarrow t,\;s_n\geq t,\ \exists x_n\in
\{\varphi(s_n,.)=0\},\\
& d(x_n, \{\varphi(t,.)=0\})\geq \tau_0,
\end{split}\]
where $d(x,\omega)=\inf_{y\in\omega} |x-y|$ is the usual distance.

Since all the $x_n$ belong to a compact set , we can
extract a converging subsequence 
$x_n\rightarrow x$. As $\varphi$ is continuous, $\varphi(t,x)=0$ or
$x\in \{\varphi(t,.)=0\}$. On the other 
hand one would also have $d(x, \{\varphi(t,.)=0\})\geq \tau_0$ which
is contradictory. 

Therefore the result of the lemma is obviously true if $t$ is such that 
$\tau_{t_0}(t)<\delta$.
 
{\em Step 2: The connection between $\tau_t$ and the Lemma.}\\
One does not have in general a uniform control on $\tau_{t_0}$. In that case, one would find a sequence $t_n$, a number $\delta>0$ s.t. $\forall \eta>0$, $\limsup \tau_{t_n}(\eta)>\delta$. 

The result of Lemma \ref{continuitylimit} would precisely rule this out and we argue by contradiction. Denote by $t$ the first time when such a jump occurs. That means that 
\begin{itemize}
\item $\tau$ can be chosen uniform till $t$, {\em i.e.} $\exists \tau$ s.t. $\forall s<t$ and $\forall \eta\in [0,\ t-s]$, $\tau_s(\eta)\leq \tau(\eta)$. 
\item There is a jump at $t$ of size $\delta>0$, which means $\exists x_0$, $\exists t_0<t$ s.t. $\varphi(t,x_0)=0$ but
\[
d(x_0,\ \{\varphi(s,.)=0\})\geq \delta,\quad\forall s\in [t_0,\ t].
\]
\end{itemize}
Note that one can take $t_0$ as close to $t$ as one wishes and in particular we may freely assume that $t-t_0$ is small enough s.t. $\tau(t-t_0)<\delta/8$.

{\em Step 3: The contradiction.}\\
First of all, note that as $\pe\rightarrow \varphi$ in $L^\infty$ norm, one has that for $\ep$ small enough and some interval $I_\ep$ with $|I_\ep|\rightarrow 0$ as $\ep\rightarrow 0$
\[
d(x_0,\ \{\pe(s,.)\geq -\fe\})\geq \frac{\delta}{2},\quad\forall s\in [t_0,\ t]\setminus I_\ep.
\]
Indeed one observes that 
\[
\{\pe(s,.)\geq -\fe\}\subset \{\varphi(s,.)\geq -\fe-\|\pe-\varphi\|_{L^\infty}\}.
\]
And by the continuity of $\varphi$
\[
\delta(\{\varphi(s,.)\geq -\fe-\|\pe-\varphi\|_{L^\infty}\},\ \{\varphi(s,.)=0\})\longrightarrow 0,
\]
for almost every $s$.

In particular we point out that one may define $I_\ep$ with $|I_\ep|\rightarrow 0$ and $\eta_\ep\rightarrow 0$ s.t.
\begin{equation}
\delta(\{\pe(s,.)\geq -\fe\},\ \{\varphi(s,.)=0\})\leq \eta_\ep,\quad\forall s\in [t_0,\ t]\setminus I_\ep.\label{closesets}
\end{equation}
Now 
define the following interval
\[
I_1(s)=[x_0-\delta/4+(s-t_0)\,V/2,\ x_0+\delta/4+(s-t_0)\,V/2].
\]
Moreover denote by $y_0$ the closest point on the left from $x_0$ in $\{\varphi(t_0,.)=0\}$. Assume for instance that $y_0<x_0$.

Next for any $\ep$ small enough, and any $s\in [t_0,\ t]$, denote 
\[
y_\ep(s)=\sup \{x<x_0,\ \pe(s,x)\geq -\fe\}.
\] 
By \eqref{closesets}, one has that $y_\ep(s)\leq y_0+\eta_\ep+\tau(s-t_0)$ for any $s\in [t_0,\ t]\setminus I_\ep$. The same result can be obtained for the closest point on the right. 

By the assumption on $x_0$ in the Lemma, this implies that $I_1$   remains at a distance larger than $\delta/4>C\,\fe$ of $\{\pe(s,.)\geq -\fe\}$.

By the properties of $D_\eps$, we deduce that  for any $s\in[s_0,\ t]$, 
\[\begin{split}
\max_{I_1(s)} r-I\star_x u_\eps(s,.)-D_\eps&\leq r(y_\ep)-I\star_x u_\ep (s,y_\ep)-D_\ep(s,y_\ep)-C\,\delta/4\\
\leq r(y_0)-I\star_x u_\ep (s,y_0)&-D_\ep(s,y_0)-C\,(\delta/4-\eta_\ep-\tau(s-t_0))\\
&\leq r(y_0)-I\star_x u_\ep (s,y_0)-D_\ep(s,y_0)-C\,\delta/16,
\end{split}\]
by choosing $\ep$ small enough. 

Now let us apply Lemma \ref{minmaxlemma} to $I_1$  to deduce that
\[\begin{split}
\max_{I_1(t)} \pe(t,.)\leq &\int_{t_0}^t (r(y_0)-I\star_x u_\ep (s,y_0)-D_\ep(s,y_0))\,ds\\
&-C\,\frac{\delta}{16}\,(t-t_0)+ C\,|I_\ep|+\max_{[x_0-3\delta/4,\ x_0+3\delta/4]} \pe(t_0,.).
\end{split}\]
Note that by \eqref{infHeps}
\[
\partial_t \pe(s,y_0)\geq r(y_0)-I\star_x u_\ep (s,y_0)-D_\ep(s,y_0)-\frac{C\ep}{\fe}.
\]
So
\[\begin{split}
\varphi(t,x_0)\leq  \|\varphi-\pe\|_{L^\infty}+ \max_{I_1(t)} \pe(t,.)\leq & \pe(t,y_0)-\pe(t_0,y_0)
-C\,\frac{\delta}{16}\,(t-t_0)\\+C\,|I_\ep|
+ \|\varphi-\pe\|_{L^\infty}+ &\max_{[x_0-3\delta/4,\ x_0+3\delta/4]} \pe(t_0,.).
\end{split}
\]
To conclude, just observe that $\varphi(t_0,y_0)=0$ so that
\[
\pe(t,y_0)-\pe(t_0,y_0)+C\,|I_\ep|\leq \|\varphi-\pe\|_{L^\infty}+\ep\,|\log \ep|+C\,|I_\ep|\longrightarrow 0,\quad as\ \ep\rightarrow 0.
\]
Therefore taking $\ep$ small enough, one concludes that
\[
\varphi(t,x_0)\leq \max_{[x_0-3\delta/4,\ x_0+3\delta/4]} \varphi(t_0,.)+\|\varphi-\pe\|_{L^\infty}<0,
\]
which gives the desired contradiction.\cqfd

%
\subsection{Control of the oscillations}
%
Before obtaining the continuity of the set $\{\pe\geq -\fe\}$ which is our main goal, we have to control the oscillations in time of the reaction term.
\begin{lemma} 
There exists a positive function $\alpha_\ep(s)$ with $\int_0^T \alpha_\ep^2(s)\,ds\leq 1$, s.t. for any $t_0\in [0,\ T]$  and any $t>t_0$
\[
\|I\star u_\ep(t,.)-I\star u_\ep(t_0,.)\|_{L^\infty}\leq \frac{1}{\sqrt \ep}\int_{t_0}^t \alpha_\ep(s)\,ds.
\]\label{oscillationreaction}
\end{lemma}
With this control, we can first show that the reaction term $r-I\star u_\ep$ essentially vanishes on the set $\{\pe\geq -\fe\}$
\begin{lemma} There exists $\tilde \ep\rightarrow 0$, $I_\ep$ composed of at most  $C\,\ep^{-3/8}$ intervals of size $\ep^{3/4}$ s.t. for any $t_0\in [0,\ T]\setminus I_\ep$, for any $x_0$ s.t.
\[
\exists t\in [t_0,\ t_0+\ep^{3/4}],\quad \pe(t,x_0)\geq -\fe
\]
then one has for any $s\in [t_0,\ t_0+\ep^{3/4}]$,
\[
r(x_0)-I\star u_\ep(s,x_0)\leq \tilde \ep.
\]\label{reprozero}
\end{lemma}
This last lemma enables us to be very precise concerning the growth rate at any point without mutations
\begin{lemma} For any $t_0$ and $t>t_0$ with $t-t_0\geq C\,\tilde \ep$, for any $x_0$ and any $\nu>C\,\sqrt{\tilde\ep}$ then
\[\begin{split}
&\mbox{either}\quad \int_{t_0}^t (r(x_0)-I\star(s,x_0)-D_\ep(s,x_0))\,ds\leq -C(t-t_0)\,(\nu^2-\tilde \ep),  \\
&\mbox{or}\quad \int_{t_0}^t D_\ep(s,x_0)\,ds\leq C(t-t_0)\,\nu.
\end{split}\]\label{growthrate}
\end{lemma}
With this control on the growth rate, it is possible to make sure that at most times the set $\{\varphi(t,.)=0\}$ is of measure $0$ or more precisely
\begin{lemma} For any $t_0\in [0,\ T]$, define the set $\Omega_{t_0}$ by $x_0\in \Omega_{t_0}$ iff
\[
x_0\in \{\varphi(t_0,.)=0\},\qquad \exists t_n\rightarrow t,\ x_n\in \{\varphi(t_n,.)=0\}\ \mbox{s.t.}\ x_n\rightarrow x_0.
\]\label{sizesupportlimit}
Then $|\Omega_t|=0$.
\end{lemma}

We now turn to the proofs of those results 

\noindent{\bf Proof of Lemma \ref{oscillationreaction}.} The proof uses one of the entropies of the system (and a different one from \cite{CJ}). Start by the following
\[\begin{split}
&\frac{d}{dt}\Big[\int u_\ep(t,x)\,I(x-y)\,u_\ep(t,y)\,dx\,dy
-2\int r\,u_\ep(t,x)\,dx\Big]\\
&=
\frac{2}{\ep}\int (r-I\star u_\ep(t,.))\,u_\ep(t,x)\,(I\star (u_\ep(t,.)-r)\,dx\\
&+\frac{2}{\ep}\int D_\ep(t,x)\,u_\ep(t,x)\,(I\star u_\ep(t,.)-r)\,dx\\
&+\int M_\ep(u_\ep)(t,x)\,(I\star u_\ep(t,.)-r)\,dx.
\end{split}\]
Note that as $D_\ep(t,x)=0$ if $\pe\geq -\fe$ then
\[
D_\ep\,u_\ep\leq C\,e^{-\fe/\ep}\leq C\,\ep.
\]
Similarly for any smooth function $\psi$, by a simple change of variable
\[
\int M_\ep(u_\ep)\psi(x)\,dx=\frac{1}{\ep} \int K(z)\,u_\ep(t,x)\,(\psi(x-\ep\,z)-\psi(x))\,dx\leq C\,\ep\,\|\nabla\psi\|_{L^\infty}.
\]
Hence one first obtains that
\begin{equation}
\frac{1}{\ep}\int_0^T\int (r-I\star u_\ep(t,.))^2\,u_\ep(t,x)\,dx\,dt\leq C.
\label{dissipationrate}
\end{equation}
Now note that
\[\begin{split}
I\star u_\ep(t,.)-I\star u_\ep(t_0,.)=&\frac{1}{\ep}\int_{t_0}^t I\star ((r-I\star u_\ep(s,.))\,u_\ep(s,.))\,ds \\
&-\frac{1}{\ep}\,\int_{t_0}^t I\star (D_\ep(s,.)\,u_\ep(s,.))\,ds\\
&+\int_{t_0}^t I\star M_\ep(u_\ep(s,.))\,ds.
\end{split}\]
With the same bounds as before, one obtains that
\[\begin{split}
\|I\star (u_\ep(t)- u_\ep(t_0))\|_{L\infty}&\leq C\,(t-t_0)+\frac{1}{\ep}\int_{t_0}^t\int |r-I\star u_\ep(s,.)|\,u_\ep(s,.)\,dx\,ds\\
\leq C\,(t-t_0)&+\frac{1}{\sqrt{\ep}}\int_{t_0}^t \left(\frac{1}{\ep}\int (r-I\star u_\ep(s,.))^2\,u_\ep(s,.)\,dx\right)^{1/2}, 
\end{split}\]
by Cauchy-Schwartz. Defining
\[
\alpha_\ep=C+\left(\frac{1}{\ep}\int (r-I\star u_\ep(s,.))^2\,u_\ep(s,.)\,dx\right)^{1/2},
\]
allows to conclude.
\cqfd

\noindent{\bf Proof of Lemma \ref{reprozero}.} First of all, consider  the intervals $I_i=[2i\,\ep^{3/4},\ 2(i+1)\,\ep^{3/4}]$ for $i=0\ldots T\,\ep^{-3/4}/2$. As $\alpha_\ep$ is uniformly bounded in $L^2$ then there are at most $C\,\ep^{-3/8}$ indices $i$ s.t.
\[
\int_{I_i} \alpha_\ep^2(s)\,ds\geq \ep^{-3/4}\left(\int_{I_i} \alpha_\ep(s)\,ds\right)^2> \ep^{3/8},
\] 
that is
\[
\frac{1}{\sqrt{\ep}}\,\int_{I_i} \alpha_\ep(s)\,ds>\ep^{1/16}.
\]
Define $I_\ep$ as the union of the intervals $I_i=[2i\,\ep^{3/4}-2\,\ep^{3/4},\ 2(i+1)\,\ep^{3/4}+2\,\ep^{3/4}]$ for such indices $i$. 

Now for any $t_0\not\in I_\ep$ and any $t\in [t_0,\ t_0+\ep^{3/4}]$, one has that
\[
\frac{1}{\ep}\,\int_{t-\ep^{3/4}}^{t+\ep^{3/4}} \alpha_\ep(s)\,ds<\ep^{1/16}.
\]
Therefore by Lemma \ref{oscillationreaction}, $I\star u_\ep(s)$ changes by at most $\ep^{1/16}$ over the interval $[t-\ep^{3/4},\ t+\ep^{3/4}]$. 

Define
\[
\tilde \ep=\ep^{1/16}+C\,\frac{\ep}{\fe}+\ep^{1/8}.
\]

Choosing any $x_0$ s.t. $\pe(t,x_0)\geq -\fe$ and any $s\in [t-\ep^{3/4},\ t+\ep^{3/4}]$, we argue by contradiction. If
\[
r(x_0)-I\star u_\ep(s,x_0)\geq \tilde\ep.
\]
Then at any other $s'\in [t-\ep^{3/4},\ t+\ep^{3/4}]$, and in particular for $s'=t_0$,
\[
r(x_0)-I\star u_\ep(s',x_0)\geq C\,\frac{\ep}{\fe}+\ep^{1/8}.
\]
Note that $D_\ep(t,x_0)=0$ and hence
\[
\partial_t \pe(t,x_0)\geq r(x_0)-I\star u_\ep(t_0,x_0)-C\,\frac{\ep}{\fe}\geq \ep^{1/8}.
\]
By continuity that means that $\pe(s',x_0)$ is increasing in $s'$ over a (possibly very small) interval $[t,\ t+\nu]$. Define $\nu$ such as to make this interval maximal. At $t+\nu$, $\pe(t+\nu,x_0)>\pe(t,x_0)\geq -\fe$ which means that $D_\ep(t+\nu,x_0)=0$ and if $t+\nu\leq t+\ep^{3/4}$
\[
\partial_t \pe(t+\nu,x_0)\geq r(x_0)-I\star u_\ep(t+\nu,x_0)-C\,\frac{\ep}{\fe}\geq \ep^{1/8}.
\]
This implies that $t+\nu>t+\ep^{3/4}$, and hence that $\pe(s',x_0)\geq -\fe$ on any $s'\in [t,\ t+\ep^{3/4}]$. Therefore
\[
\partial_t \pe(s',x_0)\geq \ep^{1/8},
\]
or finally
\[
\pe(t+\ep^{3/4},x_0)\geq -\fe+\ep^{7/8},
\]
which is impossible as the maximum of $\pe$ is of order $\ep\,|\log \ep|$. \cqfd

We turn to the next proof

\noindent{\bf Proof of Lemma \ref{growthrate}.} For any $s\in [t_0,\ t]\setminus I_\ep$, one has either that $x_0$ is at distance less than $\nu$ of $\{\pe(s,.)\geq -\fe\}$ or it is at distance larger than $\nu$. Let us decompose accordingly $[t_0,\ t]$ into $I_\ep\cap I_\nu\cap J_\nu$ where $I_\nu$ consists of the set where $x_0$ is at distance less than $\nu$ of $\{\pe(s,.)\geq -\fe\}$. 

We observe that on $J_\nu$ by the properties of $D_\ep$ and since on $\{\pe(s,.)\geq -\fe\}$, one has by Lemma \ref{reprozero} that $r-I\star u_\ep(s,.)$ is less than $\tilde \ep$ then
\[
r(x_0)-I\star u_\ep(s,x_0)-D_\ep(s,x_0)\leq \tilde\ep-C\,\nu\leq -C\,\nu.
\]
On the other hand, on $I_\nu$ then still by \ref{reprozero}, one has that
\[
r(x_0)-I\star u_\ep(s,x_0)-D_\ep(s,x_0)\leq \tilde \ep,
\]
and moreover
\[
D_\ep(s,x_0)\leq C\,\nu.
\]
Therefore
\[
\int_{t_0}^t (r(x_0)-I\star u_\ep(s,x_0)-D_\ep(s,x_0))\,ds\leq |I_\ep|+C\,\tilde \ep\, |I_\nu|-C\,\nu\,|J_\nu|,
\]
and
\[
\int_{t_0}^t D_\ep(s,x_0)\leq C\,\nu\,|I_\nu|+C\,|J_\nu|.
\]
We recall that $|I_\ep|\leq \ep^{3/8}$ which is asymptotically much smaller
than $\tilde \ep^2$. Therefore either $|J_\nu|\geq (t-t_0)\,\nu$ and we find the first possibility or $|J_\nu|\leq (t-t_0)\,\nu$ and we find the second one. \cqfd

The next proof is

\noindent{\bf Proof of Lemma \ref{sizesupportlimit}.} Note first of all that by Lemma \ref{continuitylimit}, one may in fact take any sequence $t_n\rightarrow t_0$ for any $x_0\in \Omega_{t_0}$.

Choose any $t_0$, any $\eta>0$ and any $x_0$ in $\{\varphi(t_0,.)=0\}$. Observe that for any $t>t_0$
\[
\pe(t,x_0)\geq \pe(t_0,x_0)+\int_{t_0}^t (r-I\star u_\ep(s,x_0)-D_\ep(s,x_0))\,ds-C\,\frac{\ep}{\fe}.
\]
Therefore we first deduce that
\begin{equation}
\int_{t_0}^t (r-I\star u_\ep(s,x_0)-D_\ep(s,x_0))\,ds\leq C\,\ep\,|\log \ep|+C\,(t-t_0)\,\frac{\ep}{\fe}+\|\pe-\varphi\|_{L\infty}. \label{upperbound}
\end{equation}

Next find $x_t$ the closest point from $x_0$ s.t. $\varphi(x_t)=0$. By definition $x_t-x_0\rightarrow 0$. We define the small interval around $x_t$ for $t_0\leq s\leq t$
\[
I_s=[x_t-(t-s)V/2,\ x_t+(t-s)V/2].
\]
By Lemma \ref{minmaxlemma}
\[\begin{split}
-\|\varphi-\pe\|_{L^\infty}\leq &\ep\,|\log \ep|+\int_{t_0}^t (r(x_0)-I\star_x u_\ep(s,x_0))\,ds\\
&+C\,(\bar\ep+|x_0-x_t|+(t-t_0))\,(t-t_0).
\end{split}\]
This implies that
\begin{equation}
\begin{split}
\int_{t_0}^t (r(x_0)-I\star_x u_\ep(s,x_0))\,ds\geq & -\|\varphi-\pe\|_{L^\infty}-\ep\,|\log \ep|\\
&-C\,(\bar\ep+|x_0-x_t|+(t-t_0))\,(t-t_0).
\end{split}
\label{lowerbound}
\end{equation}
Now define $\nu$ s.t. 
\[
C (\nu^2-\ep)=\frac{\|\varphi-\pe\|_{L^\infty}}{t-t_0}+\frac{\ep\,|\log \ep|}{t-t_0}+C\,(\bar\ep+|x_0-x_t|+(t-t_0)).
\]
Note that as long as $\ep$ is small enough one has that $\nu\geq \tilde \ep$. Moreover $\nu$ can be chosen arbitrarily small by taking $t$ close enough to $t_0$. By applying Lemma \ref{growthrate} for this $\nu$, \eqref{lowerbound} means that we are necessarily in the second case which means 
\[
\int_{t_0}^t D_\ep(s,x_0)\,ds\leq C\,(t-t_0)\,\nu.
\] 
We combine this inequality to \eqref{upperbound} to deduce that
\begin{equation}
\int_{t_0}^t (r-I\star u_\ep(s,x_0))\,ds\leq C\,(t-t_0)\,\left(\frac{\ep}{\fe}+\nu^2\right). \label{upperbound2}
\end{equation}
Now, we define the average
\[
\mu_\ep=\frac{1}{t-t_0}\int_{t_0}^t u_\ep(s,x_0)\,ds
\]
and from \eqref{lowerbound} and \eqref{upperbound2}, we get
\[
\left| r(x_0)-I\star\mu_\ep(x_0)\right|\leq C\,\left(\frac{\ep}{\fe}+\nu^2\right),\quad \forall x_0\in\Omega_{t_0}. 
\]
By the fundamental property \eqref{rootsnumber}, we deduce that
\[
|\Omega_{t_0}|\leq S\left(C\,\left(\frac{\ep}{\fe}+\nu^2\right)\right).
\]
To conclude take $t\rightarrow t_0$ s.t. $\nu\rightarrow 0$ and then $\ep$ small enough to obtain that $|\Omega_{t_0}|=0$. 
\cqfd

%
\subsection{Continuity in time of the set $\{x\,|\pe\geq -\fe\}$}
%
As a first consequence of our control over the oscillations of the reaction term, one may simply obtain the continuity in time of the set $\{\pe\geq -\fe\}$
\begin{lemma}
There exists $\tilde I_\ep$ with $|\tilde I_\ep|\rightarrow 0$ as $\ep\rightarrow 0$ s.t. if for some $t_0\in [0,\ T]\setminus \tilde I_\ep$, any point $x_0$, any $t>t_0$ and $\delta>0$ 
\[
d(x_0,\,\{\pe(s,.)\geq -\fe\})\geq \delta,\qquad \forall\ s\in[t_0,\ t] ,
\]
then for $\ep$ small enough with respect to $\delta$, 
\[
\pe(t,x_0)<-\fe-C\,\delta\,(t-t_0).
\]
for some constant $C$ independent of $\ep$.\label{continuity}
\end{lemma}
Unfortunately Lemma \ref{continuity} does not imply a result like \eqref{deftau} with a function $\tau$ uniform in $\ep$. It guarantees that there cannot be jumps of significant size in the set $\{\pe\geq -\fe\}$ but that set could still be propagated very fast. 

However we can combine it with Lemma \ref{sizesupportlimit} to finally deduce  a uniform (in $\ep$) continuity in time for the support of $\{\pe\geq -\fe\}$.
\begin{lemma}
 There exists $\bar I_\ep$ with $|\bar I_\ep|\rightarrow 0$ as $\ep\rightarrow 0$ s.t. for any $t_0\in [0,\ T]\setminus \bar I_\ep$, $\forall \delta>0$, there exists $\tau>0$ s.t.  $\forall x_0$, if
\[
d(x_0,\,\{\pe(t_0,.)\geq -\fe \})\geq \delta,
\]
then for $\ep$ small enough with respect to $\delta$ and $\tau$, for any $t<t_0+\tau$
\[
\pe(t,x_0)<-\fe-C\,(t-t_0)<-\fe.
\]
\label{finitespeed}
\end{lemma}

Lemma \ref{finitespeed} in particular directly implies 
\begin{corollary}
For any $t_0\in [0,\ T]\setminus \bar I_\ep$, $\exists \tau_{t_0}\in C([0,\ T])$ with $\tau_{t_0}(0)=0$ s.t. for any $t\geq t_0$ 
\[
\delta(\{\pe(t,.)\geq -\fe\},\ \{\pe(t_0,.)\geq -\fe\})\leq \tau_{t_0}(t-t_0).
\]\label{contsetep}
\end{corollary}

We turn to the proofs

\noindent{\bf Proof of Lemma \ref{continuity}.} Let first characterize the ``good'' $t_0$. We define 
\begin{equation}
\tilde I_\ep=\{t_0,\ \exists t>0\ \mbox{with}\ |I_\ep\cap [t_0,\ t]|\geq \bar\ep\,|t-t_0|\}.\label{deftildeIep}
\end{equation} 
Note that of course $I_\ep\subset \tilde I_\ep$ (just take $t=t_0+\ep^{3/4}$). 
We can estimate $|\tilde I_\ep|$ the following way: If $t_0\in \tilde I_\ep$ then there exists an index $k\geq 0$  s.t. $t_0\in \tilde I_\ep^{k}$ where $\tilde I_\ep^k$ is defined as the union of the intervals $[i\,2^k\,\ep^{3/4},\ (i+1)\,2^k\,\ep^{3/4}]$ for those $i$ s.t.
\[
 |I_\ep\cap [i\,2^k\,\ep^{3/4},\ i\,2^k\,\ep^{3/4}+2^k\,\ep^{3/4}]|\geq 2^{k-3}\,\bar\ep\,\ep^{3/4}.
\]
Each $\tilde I_\ep^k$ is composed of at most $n_k$ such intervals and since $|I_\ep|\leq C\,\ep^{3/8}$, one has that
\[
n_k\leq 2^{3-k}\,C\,\ep^{-3/8}\,\bar \ep^{-1}.
\]
Consequently
\[
|\tilde I_\ep^k|\leq C\,\ep^{3/8}\,\bar\ep^{-1}.
\]
Since $k\leq C\,|\log \ep|$, we conclude that
\[
|\tilde I_\ep|\leq C\,\ep^{3/8}\,\bar\ep^{-1}\,|\log \ep|,
\]
which indeed converges to $0$ as $\ep\rightarrow 0$.

Now define $y_\ep^{-}(s)$ similarly as in the proof of Lemma \ref{continuitylimit} by
\[
y_\eps^{-}(s)=\sup\{y<x_0,\ \pe(s,y)\geq -\fe\},
\]
and $y_\ep^{+}$ by
\[
y_\eps^{+}(s)=\sup\{y>x_0,\ \pe(s,y)\geq -\fe\}.
\]
Define again the interval
\[
I_1(s)=[x_0-\delta/2+(s-t_0)\,V/2,\ x_0+\delta/2+(s-t_0)\,V/2].
\]
It remains non empty for $s\in [t_0,\ t]$.


For $s\in [t_0,\ t]\setminus I_\ep$, one has that
\[\begin{split}
\max_{I_1(s)} (r-I\star_x u_\eps(s,.)-D_\eps)&\leq 
\max_\pm r-I\star_x u_\ep (s,.)-D_\ep(s,.)\big|_{y_\ep^\pm(s)}\\
&\qquad-C\,\delta.\\
\end{split}\]
Apply now Lemma \ref{reprozero}, if $s\not\in I_\ep$ then this implies that
\[
\max_{I_1(s)} (r-I\star_x u_\eps(s,.)-D_\eps)\leq \tilde \ep-C\,\delta\leq -C\delta,
\]
provided that $\ep$ is small enough with respect to $\delta$.

We now apply Lemma \ref{minmaxlemma} to $I_1$  to deduce that
\[\begin{split}
\max_{I_1(t)} \pe(t,.)\leq &
2\,\bar\ep(t-t_0)-C\,\delta\,(t-t_0-|I_\ep\cap[t_0,\ t]|)\\
&+C\,|I_\ep\cap [t_0,\ t]|+\max_{[x_0-3\delta/4,\ x_0+3\delta/4]} \pe(t_0,.).
\end{split}\]
Now since $t_0\not\in \tilde I_\ep$ then
\[
|I_\ep\cap [t_0,\ t]|\leq \bar \ep (t-t_0),\quad t-t_0-|I_\ep\cap[t_0,\ t]|\geq \frac{t-t_0}{2}.
\]
Hence
\[
\max_{I_1(t)} \pe(t,.)\leq -\fe-(t-t_0)(C\delta-C\,\bar\ep)\leq -\fe -C\,\delta\,(t-t_0),
\]
provided that $\ep$ is small enough with respect to $\delta$, 
which concludes the proof.\cqfd

\medskip

\noindent{\bf Proof of Lemma \ref{finitespeed}}.
By Lemma \ref{continuitylimit} and Lemma \ref{sizesupportlimit}, we know that for $a.e.\ t_0$ then
\[
|\{\varphi(t_0,.)\}|=0.
\]
Let us choose such a $t_0$ with as well $t_0\not\in \tilde I_\ep$ and denote the corresponding interval $\bar I_\ep$. Now for any $\delta$, divide the domain in intervals $I_i$ of size $\delta/2$. And denote
\[
m=\sup_i \inf_{I_i} \varphi(t_0,.).
\] 
By the previous property $m<0$. Define $\tau=-m/(3\|\varphi\|_{W^{1,\infty}})$.

Now for any $x_0$ at distance larger than $\delta$ from $\{\pe(t_0,.)\geq -\fe\}$, there exists $y_1<x_0$ and $y_2>x_0$ s.t.
\[
|y_i-x_0|\geq \delta/2,\quad \varphi(t_0,y_i)\leq m.
\]
Obviously for any $t<t_0+\tau$, 
\[
\varphi(t,y_i)\leq 2m/3,
\]
and hence for any $y$ s.t. $|y-y_i|\leq \tau$ then
\begin{equation}
\varphi(t,y)\leq m/3,\quad \pe(t,y)\leq m/3+\|\varphi-\pe\|_{L^\infty}<-\fe,
\label{propy}
\end{equation}
provided $\ep$ is chosen small enough. 

Now we apply a first time Lemma \ref{continuity} with $\delta=\tau$ which means that for $\ep$ small enough with respect to $\tau$, there cannot be a jump of size $\tau$ in $\{\pe\geq -\fe\}$. If for some $t\leq t_0+\tau$
\[
[y_1,\ y_2]\cap \{\pe(t,.)\geq -\fe\}\neq\emptyset,
\]   
then by Lemma \ref{continuity} with $\delta=\tau$, it means that for some $s\in [t_0,\ t]$, one can find a point $y$ with $|y-y_1|\leq \tau$ or $|y-y_2|\leq \tau$ s.t.
\[
\pe(s,y)\geq -\fe,
\]
which contradicts \eqref{propy}.

Hence we know that $x_0$ stays at distance $\delta/2$ of the set $\{\pe\geq -\fe\}$. We apply once more Lemma \ref{continuity} to conclude.\cqfd

\subsection{Passing to the limit in the equation: Final steps}
The first step is to characterize the limit of $D_\ep$. Of course there is a natural candidate for that which is
\begin{equation}
D=\min(K\,d(x,\ \{\varphi=0\}),\ D_0).\label{limdep}
\end{equation}
By \eqref{closesets}, we essentially know that at most times $t$, the set $\{\pe(t,.)\geq -\fe\}$ is included in a neighborhood of $\{\varphi(t,.)=0 \}$.
However to prove \eqref{limdep}, we need to show the contrary namely that $\{\varphi(t,.)=0 \}$ is included in a neighborhood of $\{\pe(t,.)\geq -\fe\}$ 
at most times $t$. This is what the following lemma shows
\begin{lemma} For any fixed $\delta>0$, denote
\[
J_\ep=\{t\in [0,\ T],\quad \delta(\{\varphi(t,.)=0\},\ \{\pe(t,.)\geq -\fe\})>\delta\}.
\]
Then $J_\ep\rightarrow 0$ as $\ep\rightarrow 0$ and consequently \eqref{limdep} holds.\label{lemmadep}
\end{lemma}
The last step is to identify the weak limit $u$ of $u_\ep$. Fortunately we now have all the required tools and we can follow the same ideas as in \cite{CJ}.
Therefore let us define
\[
 \mu_t=\mu(t)=\mu(\{\varphi(t,.)=0\}),
\] 
where $\mu$ is given by Prop.~\ref{metastable}.

We prove the intermediary result
\begin{lemma} For any fixed $t_0\in [0,\ T]\setminus \bar I_\ep$, there exists  a function
  $\sigma \in
  C(\R_+)$ with $\sigma(0)=0$ 
s.t.
\[
\int_{t_0}^t\int_\R |I\star u_\ep(s,x)-I\star\mu(t_0,x)|^2\,dx\,ds\leq
C\,(t-t_0)\,\tau(t-t_0)+\sigma(\eps),
\]
where $\tau$ is given by \eqref{deftau}.
\label{conttime}
\end{lemma}
Then simply by passing to the limit one deduces that
\begin{corollary}  For any fixed $t_0$, 
\[
\int_{t_0}^t\int_\R |I\star u(s,x)-I\star\mu(t_0,x)|^2\,dx\,ds\leq
(t-t_0)\,\tau(t-t_0).
\]
\label{conttimelimit}
\end{corollary}
This means that any $t_0$ is a Lebesgue point of $I\star u$ and then  necessarily that 
\[
I\star u(t_0,.)=I\star \mu(t_0,.),
\]
which as $\hat I>0$ lets us conclude that $u=\mu$. We have identified all the terms in the limiting equation which reads
\begin{equation}
\partial_t \varphi=r-I\star\mu(\{\varphi(t,.)=0\})-\min(K\,d(x,\ \{\varphi=0\}),\ D_0)+H(\partial_x \varphi),\label{finaleq}
\end{equation}
which finishes the proof of our main theorem. It only remains to give the proofs of Lemmas \ref{lemmadep} and \ref{conttime}.

\noindent{\bf Proof of Lemma \ref{lemmadep}.} Choose a time $t_0\in [0,\ T]\setminus \bar I_\ep$. Apply  Corollary \ref{contsetep} and Lemma \ref{finitespeed} and find $t$ s.t. $\tau_{t_0}\,(t-t_0)\leq \delta/2$ and the $\tau$ from Lemma \ref{finitespeed} satisfies $t<t_0+\tau$. 

Then for any $s\in [t_0,\ t]$ and any $x_0$ s.t.
\[
d(x_0,\ \{\pe(s,.)\geq -\fe\})>\delta.
\] 
First by  Corollary \ref{contsetep},
\[
d(x_0,\ \{\pe(t_0,.)\geq -\fe\})>\delta/2.
\] 
One has then by Lemma \ref{finitespeed} that on $[t_0,\ t]$, 
\[
\pe(s,x_0)<-\fe-C\,\delta\,(s-t_0).
\]
That implies that on $s\in [t_0+C\,\delta^{-1}\,\|\varphi-\pe\|_{L^\infty},\ t]$ then 
\[
x_0\not\in \{\varphi(s,.)=0\}.
\]
Therefore for any $s\in [t_0+C\,\delta^{-1}\,\|\varphi-\pe\|_{L^\infty},\ t]$
\[
\delta(\{\varphi(s,.)=0\},\ \{\pe(s,.)\geq -\fe\})\leq \delta.
\]
We conclude that for any $t_0\not\in \bar I_\ep$ there exists $t>0$ chosen uniformly in $\ep$ s.t.
\[
J_\ep\cap [t_0,\ t]\subset [t_0,\ t_0+C\,\delta^{-1}\,\|\varphi-\pe\|_{L^\infty}].
\]
which implies that $|J_\ep|\rightarrow 0$. \cqfd

\bigskip

\noindent {\bf Proof of Lemma \ref{conttime}.} The steps are mostly the same as in  \cite{CJ}. We give them for the sake of completeness but we do not repeat all the details here. Part of the proof can also be simplified by using the more precise results that were proved before in the present article.

\medskip

{\em Step 1: The functional.}\\ 
We look at the evolution of
\[
F_\eps(t)=\int_\R \log
u_\eps(t,x)\,d\mu_{t_0}(x)=\frac{1}{\varepsilon}\int_\R
\pe(t,x)\,d\mu_{t_0}(x), 
\]
for $s\geq t_0$.  Compute
\[\begin{split}
\frac{d}{dt} F_\eps(t)=&\frac{1}{\eps}\int_\R \left(r-I\star
u_\eps-D_\eps\right)\,d\mu_{t_0}(x)\\
&+\frac{1}{\eps}\int_\R 
H_\eps(\pe(t))\,u_\eps\,d\mu_{t_0}.
\end{split}\]
Write
\[\begin{split}
\frac{1}{\eps}\int_\R \left(r-I\star u_\eps-D_\eps
\right)&\,d\mu_{t_0}(x)=\frac{d}{dt}\int_\R
u_\eps(t,x)\,dx\\
-
\frac{1}{\eps}&\int_\R \left(r-I\star
u_\eps-D_\eps\right)\,(u_\eps(t,x)\,dx-d\mu_{t_0}(x)).  
\end{split}\]
As $r-I\star\mu_{t_0}$ vanishes on the support of
$\mu_{t_0}$, 
\[\begin{split}
  \frac{1}{\eps}\int_\R &\left(r-I\star u_\eps-D_\eps
\right)\,d\mu_{t_0}(x)=\frac{d}{dt}\int_\R 
  u_\eps(t,x)-\frac{A(t)}{\eps}\\
  &\qquad +\frac{1}{\eps}\int_\R
  I\star\left(u_\eps-\mu_{t_0}\right)\,(u_\eps(t,x)\,dx-d\mu_{t_0}(x))+\frac{1}{\ep}
\int_\R
D_\eps\,u_\eps\,dx, 
\end{split}\]
with
\[
A(t)=\int_\R \left(r-I\star\mu_s\right)\,u_\eps(t,x)\,dx+\int_\R
D_\eps d\mu_{t_0}.
\]
Notice that 
\[
\int_\R
  I\star\left(u_\eps(t,.)-\mu_{t_0}\right)\,(u_\eps(t,x)\,dx-d\mu_s(x))
\]
controls $\|I\star (u_\eps(t,.)-\mu_{t_0})\|_{L^2}^2$ since $\hat I^2\leq C\hat I$.

So we deduce since $D_\ep\geq 0$
\begin{equation}\begin{split}
\frac{1}{\eps}\int_{t_0}^t \|I\star (u_\eps(s,.)-\mu_{t_0})\|_{L^2}^2\,ds
\leq &\int_\R
\log\frac{u_\eps(t,x)}{u_\eps(t_0,x)}\,d\mu_{t_0}\\
-\int_\R 
(u_\eps(t,x)\!-\!u_\eps(t_0,x))\,dx 
&+\int_{t_0}^t\frac{A(s)}{\eps}\,ds-\frac{1}{\eps}\int_{t_0}^t\int_\R 
H_\eps(\pe(t_0))\,d\mu_{t_0}. 
\end{split}\label{conteps}\end{equation}

\medskip

{\em Step 3: Easy bounds.}\\ Lemma~\ref{apriori} tells that
\[
-H_\eps(\pe)\leq C\,\frac{\eps}{\fe}.
\]
The total mass stays bounded in time so
\[
-\int_\R 
(u_\eps(t,x)\!-\!u_\eps(t_0,x))\,dx\leq \int_\R 
(u_\eps(t,x)\!+\!u_\eps(t_0,x))\,dx\leq C.
\]
And 
\[\begin{split}
\int_\R
\log\frac{u_\eps(t,x)}{u_\eps(t_0,x)}\,d\mu_{t_0}=&\frac{1}{\eps}\int_\R
(\pe(t,x)-\pe(t_0,x))\,d\mu_{t_0}\\
&\leq\frac{1}{\eps}\int_\R
(\varphi(t,x)-\varphi(t_0,x))\,d\mu_{t_0}
+\frac{2}{\eps}\,\|\pe-\varphi\|_{L^\infty(\Omega)}.
\end{split}\] Note that by
Prop.~\ref{metastable}, $\mu_{t_0}$ is supported 
on $\{\varphi(s,.)=0\}$. Since in addition we know that $\varphi\leq
0$,
 \[\begin{split}
\int_\R
\log\frac{u_\eps(t,x)}{u_\eps(t_0,x)}\,d\mu_{t_0}
\leq \frac{2}{\eps}\,\|\pe-\varphi\|_{L^\infty(\Omega)}.
\end{split}\]
Consequently 
we deduce from \eqref{conteps} the bound
\begin{equation}
\frac{1}{\eps}\int_{t_0}^t\|I\star(u_\eps(s,.)-\mu_{t_0})\|_{L^2}^2 ds\leq
\frac{C}{\fe}+\frac{2}{\eps}\,\|\pe-\varphi\|_{L^\infty(\Omega)}+\int_{t_0}^t
\frac{A(s)}{\eps}\,ds. 
\label{conteps2}\end{equation}

\medskip

{\em Step 4: Control on $A$.}\\
The key point here is that supp$\,\mu_{t_0}=\{\varphi(t_0,.)=0\}$ is very
close to 
$\{\pe\geq -\fe\}$ as implied by Lemma \ref{lemmadep}.

First of all $D_\eps(s,.)$ is Lipschitz in $x$ and vanishes on 
$\{\pe(s,.)\geq -\fe\}$, so one has that 
\[
\int_{t_0}^t
\int_\R D_\eps(s,x)\,d\mu_{t_0}\,ds\leq
C\,\int_{t_0}^t \delta(\mbox{supp}\,\mu_{t_0},\;\{\pe(s,.)\geq -\fe\})\,ds. 
\]  
By the definition of $\mu$, its support is included in $\{\varphi=0\}$ and by Lemma \ref{continuitylimit} we finally deduce that
\begin{equation}\begin{split}
\int_{t_0}^t
\int_\R D_\eps(s,x)\,d\mu_{t_0}\,ds\leq&
C\,\tau(t-t_0)\,(t-t_0)\\
&+C\,\int_{t_0}^t \delta(\{\varphi(s,.)=0\},\;\{\pe(s,.)\geq
-\fe\})\,ds. \label{bounddeps} 
\end{split}\end{equation}
This was the only additional term with respect to \cite{CJ}. The other term is treated the same 
by decomposing
\[\begin{split}
\int_{t_0}^t\int_\R (r-I\star\mu_{t_0})\,u_\eps(s,x)\,dx
\,ds=&\int_{t_0}^t\!\!\int_\R \left(r-I\star\mu_{t_0}\right)\,
u_\eps(s,x)\ind_{\pe(s,x)\leq -\fe} 
\,dx\,ds\\
+\int_{t_0}^t &\int_\R \left(r-I\star\mu_{t_0}\right)
\,u_\eps(s,x)\ind_{\pe(s,x)\geq -\fe} 
\,dx\,ds.
\end{split}\]
For the first part, note again that by \eqref{boundeta}, there exists
$R$ s.t.
$\forall |x|>R,\ r-I\star \mu_{t_0}<0.$
 
Therefore,
we may simply dominate
\[
\int_{t_0}^t \int_\R \left( r-I\star \mu_{t_0}
\right)\,u_\eps(s,x)\ind_{\pe(s,x)\leq -\fe} 
\,dx\,ds\leq C\,(t-s)\,e^{-\fe/\eps}.
\]
Concerning the second part, just note  that $\{\pe(s,.)\geq
-\fe\}\subset \{\varphi(s,.)\geq -\alpha_\ep\}$ with $\alpha_\ep=\fe+\|\varphi-\pe\|_{L^\infty}$ and that $r-I\star \mu_{t_0}$ is non positive on $\{\varphi(t_0,.)=0\}$.
Hence 
\[\begin{split}
 \int_{t_0}^t \int_\R \left( r-I\star \mu_{t_0}
\right)\,&u_\eps(s,x)\ind_{\pe(s,x)\geq -\fe} 
\,dx\,ds\\
&\leq C\,\int_{t_0}^t \delta(\{\varphi(s,.)\geq -\alpha_\ep\},\ \{\varphi(t_0,.)=0\})\,ds.
\end{split}\] 
By Lemma \ref{continuity}, this implies that
\[\begin{split}
 \int_{t_0}^t \int_\R \left( r-I\star \mu_{t_0}
\right)\,&u_\eps(s,x)\ind_{\pe(s,x)\geq -\fe} 
\,dx\,ds\leq \tau(t-t_0)\,(t-t_0)\\
&+ C\,\int_{t_0}^t \delta(\{\varphi(s,.)\geq -\alpha_\ep\},\ \{\varphi(s,.)=0\})\,ds.
\end{split}\] 
Inequality \eqref{conteps2} now becomes
\begin{equation}\begin{split}
\int_{t_0}^t&\|I\star(u_\eps(s,.)-\mu_{t_0})\|_{L^2}^2 ds\leq
\frac{C\,\ep}{\fe}+2\,\|\pe-\varphi\|_{L^\infty(\Omega)}
+C\,(t-s)\,e^{-\fe/\eps}\\
&+C\,\tau(t-t_0)\,(t-t_0)+C\,\int_{t_0}^t \delta(\{\varphi(s,.)\geq
-\alpha_\eps\},\{\varphi(s,.)=0\})\,ds\\
&+C\,\int_{t_0}^t \delta(\{\varphi(s,.)=0\},\;\{\pe(s,.)\geq
-\fe\})\,ds .  
\end{split}\label{conteps3}\end{equation}

\medskip

{\em Conclusion.} Eq.~\eqref{conteps} indeed gives
Lemma~\ref{conttime} if one defines 
\[\begin{split}
&\tilde
\sigma(\eps)=\frac{C\eps}{\fe}+2\|\pe-\varphi\|_{L^\infty(\Omega)}+C\,T\,e^{-\fe/\eps} 
\\
&\qquad+C\,\int_0^T \delta(\{\varphi(s,.)\geq
-\alpha_\eps\},\{\varphi(s,.)=0\})\,ds\\
&+C\,\int_{0}^T \delta(\{\varphi(s,.)=0\},\;\{\pe(s,.)\geq
-\fe\})\,ds.
\end{split}\]
Thus, in order to complete the proof of 
Lemma~\ref{conttime}, we only have to check that
$\tilde{\sigma}(\varepsilon)\rightarrow 0$ when 
$\varepsilon\rightarrow 0$.

It is obvious for the first three terms. For the second note that
\[
C\,\int_0^T \delta(\{\varphi(s,.)\geq
-\alpha_\eps\},\{\varphi(s,.)=0\})\,ds\longrightarrow
0\quad\mbox{as\ }\eps\rightarrow 0,
\]
since by dominated convergence it is enough that for any $s$
\[
\delta(\{\varphi(s,.)\geq
-\alpha_\eps\},\{\varphi(s,.)=0\})\longrightarrow
0.
\]
For the last term, we use Lemma \ref{lemmadep} to show that it converges to $0$ as $\ep\rightarrow 0$. Indeed for any $\delta$ one has
\[
\int_{0}^T \delta(\{\varphi(s,.)=0\},\;\{\pe(s,.)\geq
-\fe\})\,ds\leq |J_\ep|+T\,\delta,
\]
with $|J_\ep|\rightarrow 0$.
\cqfd 
\subsection{Concluding Remarks: Finite speed of propagation of the support}
It is now easy to obtain additional properties for the limit $\varphi$, for instance the finite speed of propagation of the support. We just sketch quickly here how such properties can be proved.

First by the definition of $\mu$, one has that
\[
r-I\star \mu(\{\varphi(t,.)\})\leq 0,\quad\mbox{on}\ \{\varphi(t,.)=0\}.
\]
One then deduces immediately that for any $x$
\[
r(x)-I\star \mu(\{\varphi(t,.)\})-D(t,x)\leq 0,
\]
since $D=\min (K\,d(x,\ \{\varphi(t,.)=0\}),\ D_0)$. That means by Eq. \eqref{finaleq}
\[
\partial_t \varphi\leq H(\partial_x\varphi).
\]
Just by the following the characteristics of the Hamilton-Jacobi equation, one obtains that the set $\{\varphi=0\}$ can only propagate at finite speed $V$.

Finally let us mention that some time regularity is also known on $\mu$. For instance the set $\{\varphi=0\}$ is continuous in time by Lemma \ref{continuitylimit}.  By the stability result in \cite{JR}, one deduces that $\mu(\{\varphi(t,.)\})$ is continuous in time. As a matter of fact the finite speed of propagation of $\{\varphi=0\}$ would even imply H\"older continuity of order $1/2$ for $\mu$.

\end{document}